\documentclass[12pt]{article}
\usepackage{amsmath}
\usepackage{graphicx}
\usepackage{enumerate}
\usepackage{natbib}
\usepackage{url} % not crucial - just used below for the URL 
\usepackage{setspace}
\newcommand{\blind}{0}

\RequirePackage[OT1]{fontenc}
\RequirePackage{amsthm,amsmath}
\RequirePackage[colorlinks,citecolor=blue,urlcolor=blue]{hyperref}
\usepackage{graphicx}
\usepackage{enumerate}
\usepackage{latexsym}
\usepackage{amssymb}
\usepackage{multirow}
\usepackage{float}
\usepackage{bm}
\usepackage{paralist}
\usepackage{mathabx}
\usepackage{mathtools}
\newcommand{\cD}{\mathcal{D}}
\usepackage{natbib}
\usepackage[ruled,vlined]{algorithm2e}
\usepackage{commath}
\usepackage{subcaption}
\usepackage{soul,xcolor}
\usepackage[margin=1in]{geometry}
\usepackage{subcaption}
\usepackage[inline]{enumitem}
\usepackage[figuresright]{rotating}
\usepackage{pifont}% http://ctan.org/pkg/pifont
\usepackage{threeparttable}
{\end{quotation}}
\usepackage[normalem]{ulem}

\numberwithin{equation}{section}

\def\bbR{\mathbb{R}}

\def \mcX{\mathcal{X}}

\def \mbI{\mathbb{I}}
\def \mcR{\mathcal{R}}

\setstcolor{red}

\newcommand{\Argmin}{\mathop{\mathrm{argmin}}}

\def \mcZ{\mathcal{Z}}
\def \bfx{\mathbf{x}}

\def \bfX{\mathbf{X}}
\def \bfx{\mathbf{x}}

\def \ev{\mathbb{E}}

\def \bfv{\mathbf{v}}

\newcommand{\vo}{\vec{o}\@ifnextchar{^}{\,}{}}

\newcommand{\floor}[1]{\left\lfloor #1 \right\rfloor}
\theoremstyle{plain}

\newtheorem{theorem}{\indent Theorem}
\newtheorem*{theorem*}{\indent Theorem}

\newtheorem{Assumption}{Assumption}

\theoremstyle{definition}

\DeclareSymbolFont{largesymbolsA}{U}{txexa}{m}{n}
\DeclareMathSymbol{\varprod}{\mathop}{largesymbolsA}{16}
\usepackage{xparse}% http://ctan.org/pkg/xparse
\DeclareFontFamily{U}{mathx}{\hyphenchar\font45}
\DeclareFontShape{U}{mathx}{m}{n}{
      <5> <6> <7> <8> <9> <10>
      <10.95> <12> <14.4> <17.28> <20.74> <24.88>
      mathx10
      }{}
\DeclareSymbolFont{mathx}{U}{mathx}{m}{n}
\DeclareMathSymbol{\bigtimes}{1}{mathx}{"91}

\makeatletter
\@addtoreset{equation}{section}\makeatother

\begin{document}
\setstcolor{red}
%\title{How Well Source Data Learns Target Classifiers: a Phase Transition Phenomenon}
%\title{A Computationally Efficient Classification Algorithm in Posterior Drift Model: Phase Transition and Minimax Adaptivity}

%\runtitle{Transfer Learning Nonparametric Classification}
%	\begin{aug}
%		\author{\fnms{Ruiqi} \snm{Liu}\thanksref{t1,m1}\ead[label=e4]{ruiqliu@ttu.edu}},
%		\author{\fnms{Kexuan} \snm{Li}\thanksref{k1}\ead[label=e1]{kli@math.binghamton.edu	}},
%		\and
%		\author{\fnms{Zuofeng} \snm{Shang}\thanksref{k2}
%			\ead[label=e3]{zuofeng.shang@njit.edu}}
%		\thankstext{t1}{Corresponding author: Ruiqi Liu; email: ruiqliu@ttu.edu.}
%        \runauthor{Liu et al.}
%         \thankstext{m1}{Department of Mathematics and Statistics, Texas Tech %University, TX 79409, USA.}
%          \thankstext{k1}{Department of Mathematical Sciences, Binghamton %University, NY 13850, USA.}
%          \thankstext{k2}{Department of Mathematical Sciences, New Jersey %Institute of Technology, NJ 07102, USA.}
%	\end{aug}
	
\if0\blind
{
  \title{\bf A Computationally Efficient Classification Algorithm in Posterior Drift Model: Phase Transition and Minimax Adaptivity}
  
  \author{Ruiqi Liu\\
    Department of Mathematics and Statistics, Texas Tech University \\
    Kexuan Li\\
    Department of Mathematical Sciences, Binghamton University\\
    Zuofeng Shang\thanks{
    Zuofeng Shang gratefully acknowledge NSF grants DMS 1764280 and DMS 1821157 for supporting this work.}\\
    Department of Mathematical Sciences, New Jersey Institute of Technology}
  \maketitle
} \fi

\if1\blind
{
  \bigskip
  \bigskip
  \bigskip
  \begin{center}
    {\LARGE\bf A Computationally Efficient Classification Algorithm in Posterior Drift Model: Phase Transition and Minimax Adaptivity}
\end{center}
  \medskip
} \fi

\begin{abstract}
In massive data analysis, training and testing data often come from very different sources, and their probability distributions are not necessarily identical. 
A feature example is nonparametric classification in posterior drift model where the conditional distributions of the label given the covariates are possibly different.
In this paper, we derive minimax rate of the excess risk for nonparametric classification in posterior drift model in the setting that both training and testing data have smooth distributions, extending a recent work by \cite{cai2019transfer} who only impose smoothness condition on the distribution of testing data.
The minimax rate demonstrates a phase transition characterized by the mutual relationship between the smoothness orders of the training and testing data distributions. We also propose a computationally efficient and data-driven nearest neighbor classifier which achieves the minimax excess risk (up to a logarithm factor).  Simulation studies and a real-world application are conducted to demonstrate our approach.
\end{abstract}

\noindent%
{\it Keywords:} Transfer Learning, Domain Adaptation, Computational advantage, Adaptive Rate-Optimal
\vfill

\newpage
\section{Introduction}
Despite the significant successes of conventional classification algorithms, one of their unavoidable limitations is to assume that the source (training) data and target (testing) data are identically distributed. 
In real-world scenarios, it could be difficult and expensive to obtain source data that has the same distribution as the target data (\citealp{weiss2016survey}). Thus, an algorithm which can overcome such discrepancy would be highly valuable. Transfer learning is a promising tool to build models for the \textit{target domain} by transferring data information from the related \textit{source domain}.  %A commonly used term in the literature to describe the process of transferring information between domains is called \textit{domain adaptation}.
In comparison with traditional machine learning algorithms, transfer learning has demonstrated advantages in many aspects such as  image classification (\citealp{zhu2011heterogeneous, han2018new,hussain2018study}), autonomous driving (\citealp{Kim_2017_CVPR_Workshops}), recommendation system (\citealp{zhao2013active,zhang2017cross}), etc.
An important transfer learning technique is the so-called \textit{domain adaptation} (\citealp{weiss2016survey}), in which the information or knowledge is adapted from one or more source domains to the target domain. Empirical successes of domain adaptation have attracted increasing attention to the study of its theoretical properties. For example, \cite{bendavid2007domain} derive a bound on the generalization error of the classifiers trained from data in the source domain, later extended by \cite{john2008learningbounds}, \cite{zhang2012generalization} and \cite{zhao2018adversarial} to the case where the classifiers are trained
from both source and target domains. Researchers have also
proposed additional structural relationships between the source and target domains, e.g., the \textit{covariate shift} model with different marginal distributions and the \textit{posterior drift} model with different conditional distributions, under which the generalization error bounds are successfully established. To name a few, see  
\cite{shimodaira2000improving, huang2007correcting, sugiyama2008direct, mansour2012multiple, hoffman2018algorithms, kpotufe2018marginal, scott2013classification,natarajan2013learning,manwani2013noise,gao2016risk, natarajan2017cost, cannings2020classification, cai2019transfer}.

A notable work is \cite{cai2019transfer} who propose a $k$-nearest neighbor ($k$NN) classifier
based on the posterior drift model and derive minimax optimality. 
The $k$-nearest neighbors are detected over the entire source and target data
which could be computationally expensive when data size is large. 
Meanwhile, the theoretical results 
only involve the smoothness of the target distribution but the impact of the smoothness of the source distributions remains unknown.
The aim of this work is to further strengthen the two aspects.  
Specifically, we propose a more computationally efficient $k$NN classifier that
requires detecting the nearest neighbors for local data only. We discover a phase transition phenomenon for the minimax excess risk characterized by the mutual relationship between the smoothness orders of the source and target distributions, which degenerate to \cite{cai2019transfer} in the special case when the smoothness order of the source distributions vanishes. 
Such a phenomenon provides a more complete understanding on the impact of smooth data distributions in transfer learning.
In the following subsection, we describe our contributions in more details. Before that, let us introduce some terminologies and notation. 

\textbf{Terminologies and Notation:} Let $\|\bfv\|^2=\bfv^\top\bfv$  denote the Euclidean norm of the vector $\bfv$. For two sequences $a_n$ and $b_n$, we say $a_n\lesssim b_n$ if $a_n\leq cb_n$ for some constant $c>0$ and all sufficiently large $n$. For $a>0$, let $\floor{a}$ be the largest integer that less than or equal to $a$. Let $\cD=[0,1]^d\times\{0,1\}$ and $P,Q$ be probability distributions over $\cD$.
 For $(\bfX,Y)\in\cD$, the $d$-dimensional $\bfX$ is regarded as covariates or features, and $Y$ is the binary label of $\bfX$. For $(\bfX,Y)$ drawn from $P$ (or $Q$), let $P_\bfX$ (or $Q_\bfX$) denote the marginal probability distribution of $\bfX$. Define $\eta_P(\bfx)=P(Y=1|\bfX=\bfx)$ 
 and $\eta_Q(\bfx)=Q(Y=1|\bfX=\bfx)$ as the conditional distributions of $Y$ given $\bfX=\bfx$ for any $\bfx\in[0,1]^d$. For positive sequences $a_n$ and $b_n$, we say $a_n\lesssim b_n$ ($a_n\gtrsim b_n$) if $a_n\leq cb_n$ ($a_n\geq cb_n$) for some $c>0$ and all large enough $n$. We say $a_n\asymp b_n$ if $a_n\lesssim b_n$ and $a_n\gtrsim b_n$. We use $\textrm{supp}(\mu)$ to denote the support of the a probability measure $\mu$.
  
\subsection{Nonparametric classification in posterior drift model}
Suppose that
 $(\bfX_1^P,Y_1^P),\ldots,(\bfX_{n_P}^P,Y_{n_P}^P)$
 are $n_P$ i.i.d. observations from $P$, and $(\bfX_1^Q,Y_1^Q),\ldots,$ $(\bfX_{n_Q}^Q,Y_{n_Q}^Q)$ are $n_Q$ 
  i.i.d. observations
from $Q$. Moreover, we assume these observations from $P$ and $Q$ are mutually independent.
For simplicity, we call $\{(\bfX_i^P,Y_i^P)\}_{i=1}^{n_P}$ the $P$-data
and $\{(\bfX_i^Q,Y_i^Q)\}_{i=1}^{n_Q}$ the $Q$-data.
Given future covariates $\bfX$ from $Q_\bfX$, 
we are interested in predicting its unknown binary label $Y$ based on the full training data $\mcZ$, where
\begin{equation}\label{eqn:entire:data}
\mcZ=\{(\bfX_1^P,Y_1^P),\ldots,(\bfX_{n_P}^P,Y_{n_P}^P), (\bfX_1^Q,Y_1^Q),\ldots,(\bfX_{n_Q}^Q,Y_{n_Q}^Q)
\}.
\end{equation}
We adopt the posterior drift model proposed by \cite{scott2019generalized}, in which $P_\bfX$ and $Q_\bfX$ have common supports, whereas the conditional distributions of $Y$ given $\bfX$ under $P,Q$, namely $\eta_P$ and $\eta_Q$, are possibly different. In this model, $P$ is the \textit{source distribution}
and $Q$ is the \textit{target distribution}.
This model has been recently adopted by \cite{cai2019transfer} 
in nonparametric classification who proposed an optimal adaptive $k$NN classifier. 
Since their method requires identifying the nearest covariates among all covariates in $\mcZ$, which requires $n_P+n_Q$ attempts and might be computationally expensive. The computational cost easily scales up when the training data consists of observations from multiple distributions. Hence, it is interesting to design a more efficient algorithm that can achieve the same optimality. 
\subsection{Our Contributions}
\begin{figure}
\centering
\includegraphics[width=3 in, height=2.5 in]{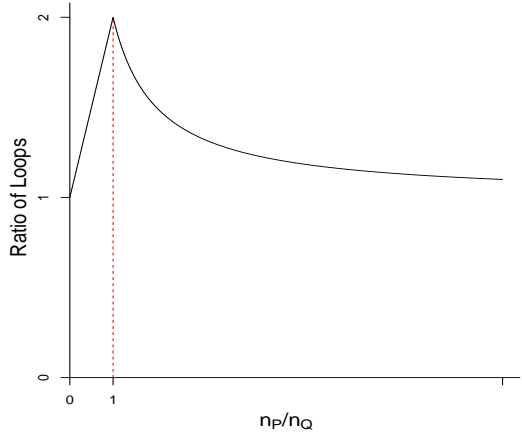}
\caption{\it The ratio of attempts required by \cite{cai2019transfer} and Algorithm \ref{alg:ag2}. }
\label{figure:runtime}
\end{figure}

Our first contribution is to propose a more computationally efficient adaptive $k$-NN classifier, i.e., Algorithm \ref{alg:ag2} in Section \ref{sec:scalable:alg}. Notably, our method requires $\max\{n_P,n_Q\}$ attempts to identify the nearest covariates. Consequently, the ratio of attempts required by \cite{cai2019transfer} and Algorithm \ref{alg:ag2} is 
 $\frac{n_P+n_Q}{\max\{n_P,n_Q\}}$, which is nearly 2 if $n_P\approx n_Q$
 (see Figure \ref{figure:runtime}).
 In other words, when $P$-data and $Q$-data have equal amount of data points,
 Algorithm \ref{alg:ag2} only requires nearly half computational cost
 of \cite{cai2019transfer}. 
The ratio of attempts further increases 
if more source distributions are involved,
hence, the computational cost of our method, compared to \cite{cai2019transfer}, will be further reduced.
See discussions in Section \ref{section:multiple:source:model}.

\begin{figure}[H]
    \centering
    \includegraphics[width=2.7in]{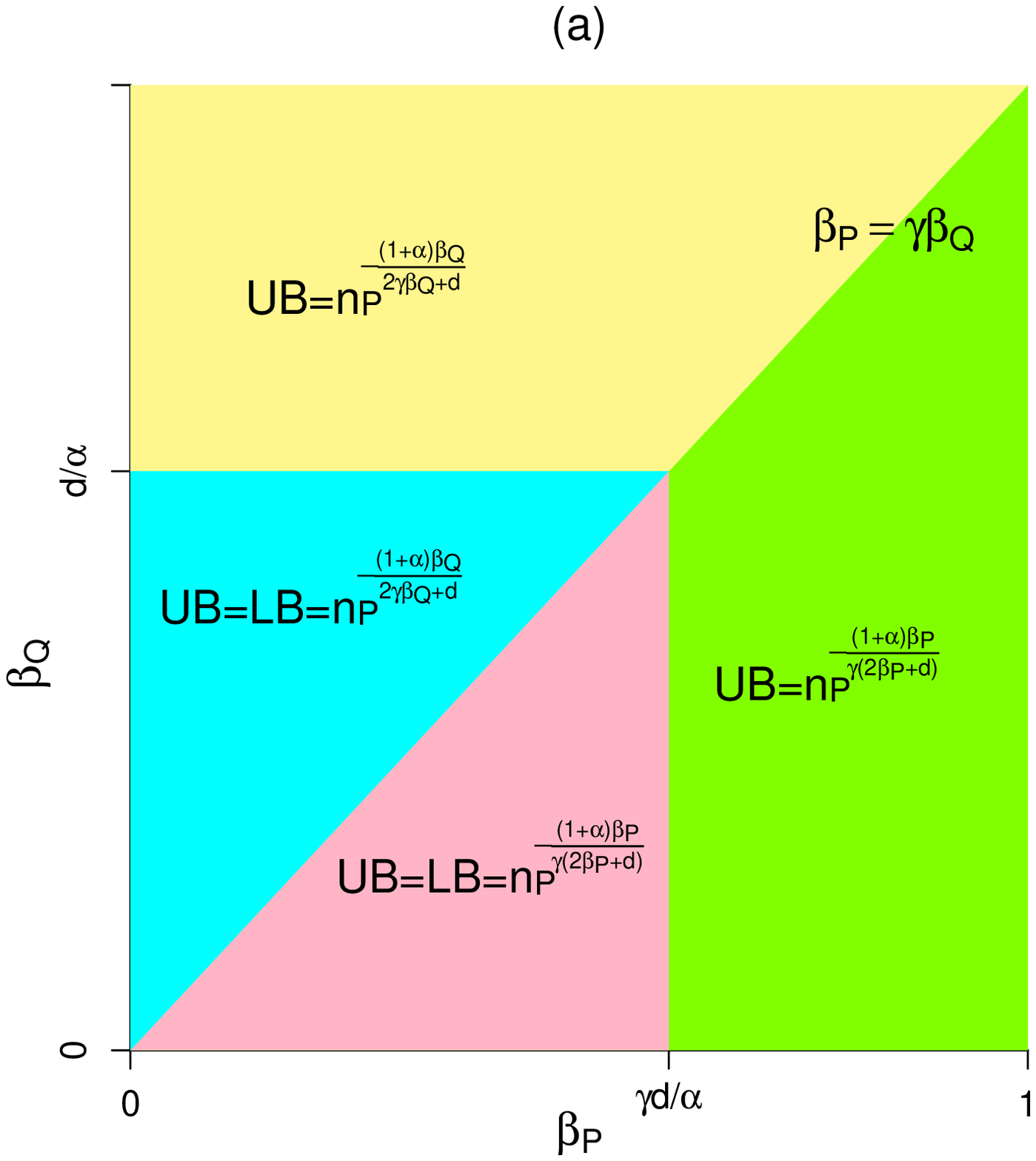}
    \hspace{10mm}
    \includegraphics[width=2.7in]{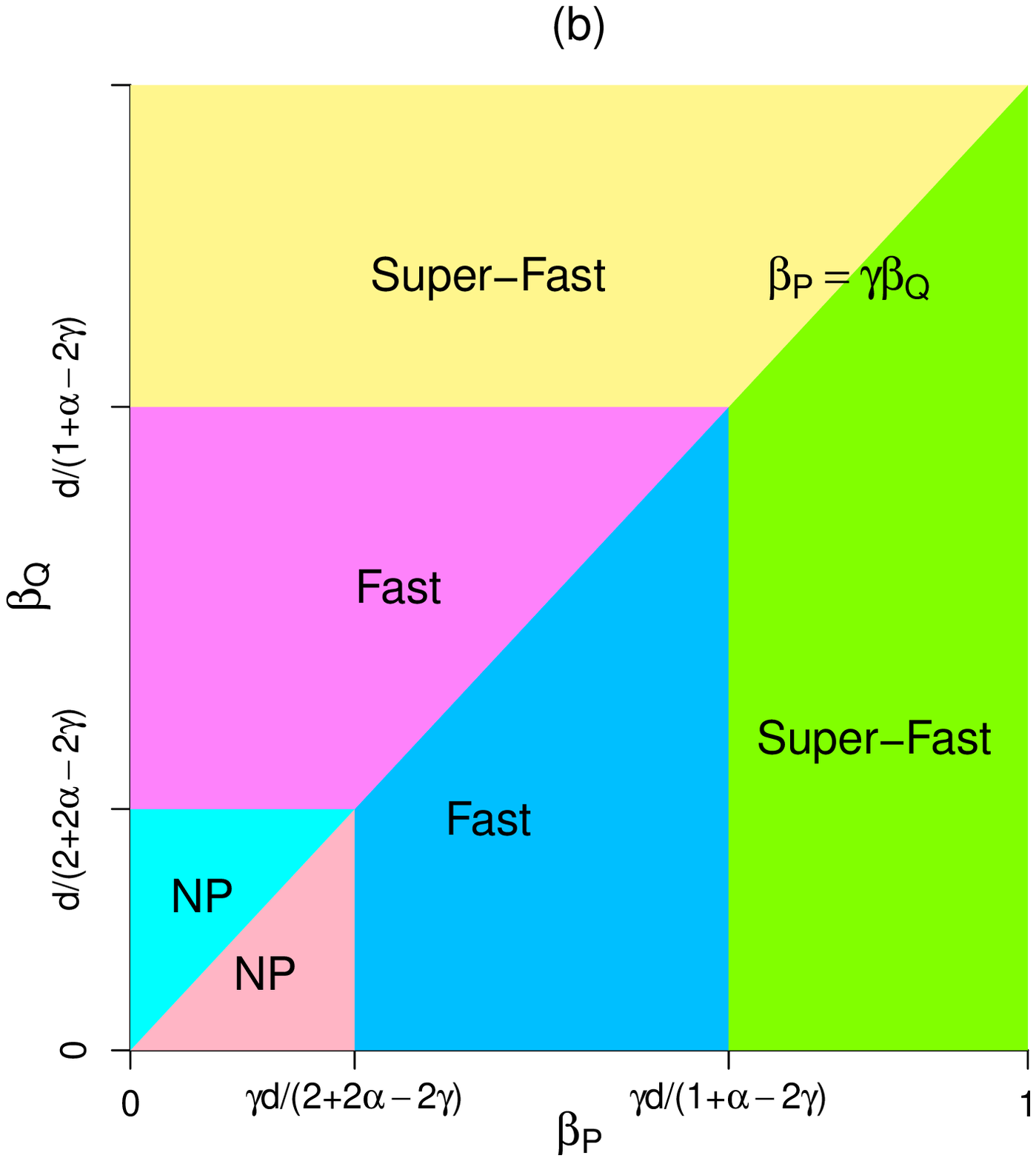}
    \caption{\it (a). Upper bounds (UB) and lower bounds (LB) of the minimax excess risk in different domains of $(\beta_P,\beta_Q)$.
(b). Categories of the excess risk bounds characterized by $(\beta_P,\beta_Q)$. NP: nonparametric rate $(\gtrsim n_P^{-1/2})$; Fast: fast rate $(\lesssim n_P^{-1/2})$; Super-Fast: super-fast rate $(\lesssim n_P^{-1})$.
}
\label{fig:phase_transition}
\end{figure}

Our second contribution is in theoretical aspect. 
We establish exact orders for the minimax excess risk when $(\eta_P,\eta_Q)$ are $(\beta_P,\beta_Q)$-H\"{o}lder smooth (see Assumption \ref{A1:holder:smooth}). In contrast, \cite{cai2019transfer} established
exact orders for the minimax excess risk in the special case $\beta_P=0$.
It turns out that when $\beta_P>0$, the minimax rate has a faster order,
demonstrating the advantage of utilizing smooth source distributions.
To describe our findings, consider $n_Q=0$ for simplicity.
Below is a summary of the results in which the orders for the minimax excess risk demonstrate a phase transition characterized by a mutual relationship between $\beta_P$ and $\beta_Q$:
\begin{equation}\label{order:motivation}
\textrm{the minimax excess risk}\asymp
\left\{
\begin{array}{cc}
n_P^{-\frac{(1+\alpha)\beta_Q}{\gamma(2\beta_P+d)}}, & 
\textrm{if $\gamma\beta_Q\le\beta_P\le \gamma d/\alpha$,}\\
n_P^{-\frac{(1+\alpha)\beta_Q}{2\gamma\beta_Q+d}},  & \textrm{if $\beta_P<\gamma\beta_Q\le \gamma d/\alpha$,} 
\end{array}
\right.
\end{equation}
where $\alpha>0$ quantifies the Tsybakov noise level (see Assumption \ref{A1:marginal:assumption}) and $\gamma>0$ measures the relative signal strength of $P$ and $Q$ (see Assumption \ref{A1:relative:signal:multipe:source}).
See Figure \ref{fig:phase_transition}(a) for an illustration of (\ref{order:motivation}).
Note that (\ref{order:motivation}) 
excludes the regions $\beta_P>\gamma d/\alpha$
and $\beta_Q>d/\alpha$ in which only upper bounds on the minimax excess risk
are available. Interestingly, the upper bounds are super fast $(\lesssim n_P^{-1})$, see Figure \ref{fig:phase_transition}(b), which is
consistent to the findings of \cite{audibert2007fast}.
The orders of the minimax excess risk are fast $(\lesssim n_P^{-1/2})$
or in a nonparametric rate ($\gtrsim n_P^{-1/2}$) in other domains of $(\beta_P,\beta_Q)$.
The results are further extended to general $n_Q$ in Section \ref{secton:two:source:model}, and to multiple source distributions 
in Section \ref{section:multiple:source:model}.
%\begin{figure}[H]
%\centering
%\begin{minipage}[b]{0.45\linewidth}
%\includegraphics[width=2.7in]{phase2.eps}
%\caption{\it Upper bounds (UB) and lower bounds (LB) of the minimax excess %risk
%in different domains of $(\beta_P,\beta_Q)$.
%Here, $\alpha$ is the Tysbakov noise exponent (see Assumption \ref{A1:marginal:assumption}) and $\gamma$ is the relative signal exponent
%(see Assumption \ref{A1:relative:signal}).}
%\label{fig:phase2}
%\end{minipage}%
%\quad\quad\quad\quad
%\begin{minipage}[b]{0.45\linewidth}
%\includegraphics[width=2.7in]{phase1.eps}
%\caption{Categories of the excess risk bounds characterized by $(\beta_P,\beta_Q)$. NP: nonparametric rate $(\gtrsim n^{-1/2})$; Fast: fast rate $(\lesssim n^{-1/2})$; Super-Fast: super-fast rate $(\lesssim n^{-1})$.}
%\label{fig:phase1}
%\end{minipage}%
%\end{figure}

\section{A Computationally Efficient Adaptive $k$NN Classifier}\label{sec:scalable:alg}
In this section, we propose a computationally efficient adaptive $k$NN classifier. 
For any $\bfx \in [0, 1]^d$, $1\le k_P\le n_P$ and $1\le i\le k_P$,
let $\bfX_{(i)}^P(\bfx)$
denote the $i$th nearest covariate of $\bfx$ among
$\bfX_1^P, \bfX_2^P,\ldots, \bfX_{n_P}^P$,
and let 
$Y_{(i)}^P(\bfx)$ denote the label of $\bfX_{(i)}^P(\bfx)$. 
One can similarly define $\bfX_{(i)}^Q(\bfx)$ and $Y_{(i)}^Q(\bfx)$.
Let
\begin{equation}\label{eq:knn:kp}
\widehat{\eta}_{k_P}(\bfx)=\frac{1}{k_P}\sum_{i=1}^{k_P}Y_{(i)}^P(\bfx),
\,\,\,\,\,\,\,\,
\widehat{\eta}_{k_Q}(\bfx)=\frac{1}{k_Q}\sum_{i=1}^{k_Q}Y_{(i)}^Q(\bfx).
\end{equation}
Here $\widehat{\eta}_{k_P}(\bfx)$ ($\widehat{\eta}_{k_Q}(\bfx)$) is the 
$k$NN estimator of $\eta_P(\bfx)$ ($\eta_Q(\bfx)$) based on the $P$-data
 ($Q$-data).
Inspired from \cite{cai2019transfer}, we can aggregate $\widehat{\eta}_{k_P}$
and $\widehat{\eta}_{k_Q}$ into a weighted $k$NN estimator:
\begin{eqnarray*}
\widehat{\eta}_{NN}(\bfx)=\frac{w_Pk_P\widehat{\eta}_{k_P}(\bfx)+w_Qk_Q\widehat{\eta}_{k_Q}(\bfx)}{w_Pk_P+w_Qk_Q},
\end{eqnarray*}
where $w_P,w_Q$ are positive weights.
The corresponding $k$NN classifier $\widehat{f}_{NN}$ is then defined as 
\begin{equation}\label{eq:knn:classifier}
\widehat{f}_{NN}(\bfx):=\mbI(\widehat{\eta}_{NN}(\bfx)\geq 1/2)=\begin{cases}
1, & \textrm{ if } \widehat{\eta}_{NN}(\bfx)\geq 1/2,\\
0, & \textrm{ if } \widehat{\eta}_{NN}(\bfx)< 1/2,
\end{cases}
\,\,\,\,\,\,\,\,
\textrm{for any $\bfx\in[0,1]^d$.}
\end{equation}
A limitation of $\widehat{f}_{NN}$ is that it requires predetermined $k_P,k_Q,w_P,w_Q$. To address this, we propose 
Algorithm \ref{alg:ag2} in which the parameters are data-driven.

\begin{algorithm}[htp!]
\SetAlgoLined
\KwIn{$P$-data $(\bfX_1^P, Y_1^P),\ldots, (\bfX_{n_P}^P, Y_{n_P}^P)$, $Q$-data $(\bfX_1^Q, Y_1^Q),\ldots, (\bfX_{n_Q}^Q, Y_{n_Q}^Q)$ with $n_P\geq n_Q$, and new features $\bfx$;}
% \textbf{Initialization:} $k=1$\;
\textbf{Initiation: } set $k_P=0$;

\textbf{while} $k_P< n_P$ \textbf{do}

\hspace*{0.5cm}  update $k_P:=k_P+1$ and $k_Q:=\floor{k_Pn_Q/n_P}$;

\hspace*{0.5cm} calculate $\widehat{\eta}_{k_P}$ and $\widehat{\eta}_{k_Q}(\bfx)$ (set $\widehat{\eta}_{k_Q}(\bfx)=1/2$ if $k_Q=0$);

\hspace*{0.5cm} calculate

\hspace*{0.5cm}   $r_{k_P}=\begin{cases}
\sqrt{k_P(\widehat{\eta}_{k_P}(\bfx)-1/2)^2+k_Q(\widehat{\eta}_{k_Q}(\bfx)-1/2)^2} & \textrm{ if } sign(\widehat{\eta}_{k_P}(\bfX)-1/2)=sign(\widehat{\eta}_{k_Q}(\bfX)-1/2)\nonumber;\\
\max\{\sqrt{k_P}|\widehat{\eta}_{k_P}(\bfx)-1/2|,\sqrt{k_Q}|\widehat{\eta}_{k_Q}(\bfx)-1/2|\}& \textrm{ if } sign(\widehat{\eta}_{k_P}(\bfx)-1/2)\neq sign(\widehat{\eta}_{k_Q}(\bfx)-1/2)\nonumber;\\
\end{cases}$

\hspace*{0.5cm} \textbf{if } $r_{k_P}>\sqrt{[d+\log(n_P+n_q)]\log(n_P+n_q)}$  \textbf{or} $k_P=n_P$ \textbf{ then} 

\hspace*{1cm}  set $\widehat{k}_P=k_P$ and $\widehat{k}_Q=k_Q$;

\hspace*{1cm}  exit  loop;

\hspace*{0.5cm} \textbf{end if}

\textbf{end while}

calculate  $\widehat{\eta}_{\widehat{k}_P}(\bfx)$, $\widehat{\eta}_{\widehat{k}_Q}(\bfx)$;

\KwOut{classifier $\widehat{f}_{pa}(\bfx)=\mbI({k_P}(\widehat{\eta}_{\widehat{k}_P}(\bfx)-1/2)+{k_Q}(\widehat{\eta}_{\widehat{k}_Q}(\bfx)-1/2)\geq 0)$.}
 \caption{An Adaptive $k$NN Algorithm}
 \label{alg:ag2}
\end{algorithm}
To ease presentation, Algorithm \ref{alg:ag2} has only considered $n_P\geq n_Q$. When $n_P<n_Q$, by flipping $n_P$ and $n_Q$ we can set $k_P=\floor{k_Qn_P/n_Q}$ during the loops until the same stopping rule is met. Below we discuss the intuition why Algorithm \ref{alg:ag2} performs optimal. Let $\rho=\frac{w_Pk_P}{w_Pk_P+w_Qk_Q}$, so we can rewrite $\widehat{\eta}_{NN}$ as $\widehat{\eta}_{NN}(\bfx)=\rho \widehat{\eta}_{k_P}(\bfx)+(1-\rho)\widehat{\eta}_{k_Q}(\bfx)$, whose ``signal,"
defined as the absolute deviation from random guess, and ``standard deviation" are given as $|\widehat{\eta}_{NN}(\bfx)-1/2|$ and $\sqrt{\rho^2/k_P+(1-\rho)^2/k_Q}$, respectively. During each loop, Algorithm \ref{alg:ag2} finds the ``optimal" $\rho$ that minimizes the ``signal-to-noise" ratio:
\begin{eqnarray}
\widehat{\rho}=\Argmin_{\rho}\frac{|\widehat{\eta}_{NN}(\bfx)-1/2|}{\sqrt{\rho^2/k_P+(1-\rho)^2/k_Q}}.\label{eq:optimal:rho}
\end{eqnarray}
By direct calculations, the minimal value of (\ref{eq:optimal:rho}) is $r_{k_P}$ which is achieved at $\rho=\widehat{\rho}$. Therefore, Algorithm \ref{alg:ag2} scans the first $k_P$ such that the minimal ``signal-to-noise" ratio $r_{k_P}$ is greater than a threshold $\sqrt{[d+\log(n_P+n_q)]\log(n_P+n_q)}$. The choice of such threshold is inspired from \cite{cai2019transfer}, under which it can be shown that, with high probability, $\ev(\widehat{\eta}_{NN}(\bfx)|\mcX)$ and $\widehat{\eta}_{NN}(\bfx)$ have the same sign. Here $\mcX:=\{\bfX_1^P,\ldots, \bfX_{n_P}^P, \bfX_1^Q, \ldots, \bfX_{n_Q}^Q\}$ is the collection of all covariates.
Under assumptions in Section \ref{section:asymptotic:theory}, is can be shown that $\ev(\widehat{\eta}_{NN}(\bfx)|\mcX)$ has the same sign as $\eta_Q(\bfx)$,
so $\widehat{\eta}_{NN}(\bfx)$ is asymptotically optimal.

To conclude this section, we briefly discuss the computing advantage of our method.
In Algorithm \ref{alg:ag2}, $(k_P,k_Q)$ is selected over $\{(k_P,k_Q): k_Q=\floor{k_Pn_Q/n_P}, 1\le k_P\le \max\{n_P,n_Q\}\}$
which requires $\max\{n_P,n_Q\}$ attempts.
In contrast, \cite{cai2019transfer} selects $(k_P,k_Q)$
over a set of cardinality $n_P+n_Q$, hence, requires $n_P+n_Q$ attempts.  
Therefore, Algorithm \ref{alg:ag2} is computationally more feasible than \cite{cai2019transfer}.

\section{Asymptotic Theory}\label{section:asymptotic:theory}

In this section, we explore the asymptotic properties of $\widehat{f}_{NN}$ in (\ref{eq:knn:classifier})
and $\widehat{f}_{pa}$ provided in Algorithm \ref{alg:ag2}. We start from the easier case $n_Q=0$, and proceed to the general case $n_Q\neq 0$.

\subsection{Minimax Rate for the Excess Risk when $n_Q=0$}\label{section:one:source}
When $n_Q=0$, we have
$\widehat{f}_{NN}(\bfx)=\mbI(\widehat{\eta}_{k_P}(\bfx)\geq 1/2)$,
where $\widehat{\eta}_{k_P}(\bfx)$ is given in (\ref{eq:knn:kp}). Meanwhile,
by (\ref{eq:optimal:rho}) we have $\widehat{\rho}=1$. Thus, Algorithm \ref{alg:ag2} becomes the following Algorithm \ref{alg:ag1}.
\begin{algorithm}[htp!]
\SetAlgoLined
\KwIn{$P$-data $(\bfX_1^P, Y_1^P),\ldots, (\bfX_{n_P}^P, Y_{n_P}^P)$ and new features $\bfx$;}
% \textbf{Initialization:} $k=1$\;
\textbf{Initiation: } set $k_P=0$;

\textbf{while} $k_P< n_P$ \textbf{do}

\hspace*{0.5cm}  update $k_P:=k_P+1$;

\hspace*{0.5cm} calculate $\widehat{\eta}_{k_P}(\bfx)=\frac{1}{k_P}\sum_{i=1}^{k_P} Y_{(i)}^P(\bfx)$;

\hspace*{0.5cm} calculate $r_{k_P}=\sqrt{k_P}|\widehat{\eta}_{k_P}(\bfx)-1/2|$;

\hspace*{0.5cm} \textbf{if } $r_{k_P}>\sqrt{[d+\log(n_P)]\log(n_P)}$  \textbf{or} $k=n_P$\textbf{ then} 

\hspace*{1cm}  set $\widehat{k}_P=k_P$;

\hspace*{1cm}  exit  loop;

\hspace*{0.5cm} \textbf{end if}

\textbf{end while}

calculate $\widehat{\eta}_{\widehat{k}_P}(\bfx)$;

\KwOut{classifier $\widehat{f}_{pa}(\bfx)=\mbI(\widehat{\eta}_{\widehat{k}_P}(\bfx)\geq 1/2)$.}
 \caption{An Adaptive $k$NN Algorithm under $n_Q=0$}
 \label{alg:ag1}
\end{algorithm}

Before investigating the asymptotic properties of the above classifiers, we introduce some  technical assumptions. Throughout, let $\lambda$ denote the 
Lebesgue measure on $\bbR^d$.
\begin{Assumption}\label{A1:strong:density} (Common Support and Strong Density)  There exist an $\Omega \subset [0, 1]^d$ and  constants $c_\lambda, r_\lambda>0$ such that
\begin{enumerate}[label=(\alph*)]
\item\label{A1:strong:density:a} $\Omega$ is the common support of the marginal distributions $P_{\bfX}$ and $Q_{\bfX}$;
\item \label{A1:strong:density:b} $\lambda [\Omega \cap B(\bfx, r)]\geq c_\lambda \lambda[B(\bfx,r)]$ for all $0<r<r_\lambda$ and $\bfx\in \Omega$;
\item \label{A1:strong:density:c}  $c_\lambda<\frac{dP_{\bfX}}{d\lambda}(\bfx)<c_\lambda^{-1}$ and $c_\lambda<\frac{dQ_{\bfX}}{d\lambda}(\bfx)<c_\lambda^{-1}$, for all $\bfx \in \Omega$.
\end{enumerate}

\end{Assumption}
\begin{Assumption}\label{A1:holder:smooth} (H\"{o}lder Smoothness) 
There exist constants  $C_\beta>0$ and $\beta_P, \beta_Q\in [0,1]$ with $\max(\beta_P, \beta_Q)>0$ such that 
$(\eta_P,\eta_Q)$ are $(\beta_P,\beta_Q)$-H\"{o}lder smooth, i.e., $|\eta_P(\bfx_1)-\eta_P(\bfx_2)|\leq C_\beta \|\bfx_1-\bfx_2\|^{\beta_P}$ and  $|\eta_Q(\bfx_1)-\eta_Q(\bfx_2)|\leq C_\beta \|\bfx_1-\bfx_2\|^{\beta_Q}$ for all $\bfx_1, \bfx_2 \in \Omega$.
%smoothness $\eta_P$ is $(\beta_P, C_\beta)$-H\"older smooth
\end{Assumption}
\begin{Assumption}\label{A1:marginal:assumption}(Tsybakov’s Noise Condition) There exist constants $\alpha\geq 0$ and $C_\alpha>0$ such that, for all $t\in (0, 1/2]$, $Q_\bfX(|\eta_Q(\bfX)-1/2|<t)\leq C_\alpha t^{\alpha}$.
\end{Assumption}
\begin{Assumption}\label{A1:relative:signal}(Relative Signal Exponent Condition) For all $\bfx \in \Omega$, it holds that
\begin{enumerate}[label=(\alph*)]
\item \label{A1:relative:signal:a} $(\eta_P(\bfx)-1/2)(\eta_Q(\bfx)-1/2)\geq 0$;
\item \label{A1:relative:signal:b} $|\eta_P(\bfx)-1/2|\geq C_\gamma|\eta_Q(\bfx)-1/2|^\gamma$ for some constants $\gamma, C_\gamma>0$.
\end{enumerate}
\end{Assumption}

Assumption \ref{A1:strong:density} consists of three aspects. Assumption  \ref{A1:strong:density}\ref{A1:strong:density:a} requires \textrm{supp}($Q_\bfX$)$=$\textrm{supp}($P_\bfX$) which can be relaxed to \textrm{supp}($Q_\bfX$)$\subset$ \textrm{supp}($P_\bfX$) with a slight modification in the proof. The latter is necessary since otherwise there will be a covariate in \textrm{supp}($Q_\bfX$)$\backslash$ \textrm{supp}($P_\bfX$) which is unpredictable by the $P$-data.
Assumptions \ref{A1:strong:density}\ref{A1:strong:density:b} and \ref{A1:strong:density}\ref{A1:strong:density:c} are the so-called \textit{Strong Density Assumption} commonly used in literature
(see \citealp{audibert2007fast}).
Assumption \ref{A1:strong:density}\ref{A1:strong:density:b} regularizes the feature space $\Omega$. Assumption \ref{A1:strong:density}\ref{A1:strong:density:c} assumes that the densities of $P_\bfX$ and $Q_\bfX$ are bounded away from zero and infinity. 

Assumption \ref{A1:holder:smooth} requires that $\eta_P$ and $\eta_Q$ are H\"{o}lder smooth with orders $\beta_P$ and $\beta_Q$, which includes  the special case $\beta_P=0$ in  \cite{cai2019transfer}. 
Assumption \ref{A1:marginal:assumption} is the so-called Tsybakov noise condition with a noise exponent $\alpha$ (see \citealp{mammen1999smooth, audibert2007fast}). 
Assumption \ref{A1:relative:signal}, firstly introduced by \cite{cai2019transfer}, consists of two parts. Assumption \ref{A1:relative:signal}\ref{A1:relative:signal:a} requires that the Bayes classifiers $f_P^*(\cdot)=\mbI(\eta_P(\cdot)\geq 1/2)$ and $f_Q^*(\cdot)=\mbI(\eta_Q(\cdot)\geq 1/2)$ are essentially the same. 
Assumption \ref{A1:relative:signal}\ref{A1:relative:signal:b} measures the relative signal strength of $P$ and $Q$. Similar assumption is also proposed by \cite{hanneke2019value}.

Recall that the excess risk of $\widehat{f}_{NN}$ under $Q$ is defined as
\begin{eqnarray}
\mcR_Q(\widehat{f}_{NN})=\ev(\mathcal{E}_Q(\widehat{f}_{NN}))-\mathcal{E}_Q(f^*_Q),\nonumber
\end{eqnarray}
where $\mathcal{E}_Q(f)=Q(f(\bfX)\neq Y)$ is the classification risk of classifier $f: \Omega \to \{0, 1\}$ under $Q$, and  $f_Q^*$ is the Bayes classifier defined as
$f_Q^*(\cdot)=\mbI(\eta_Q(\cdot)\geq 1/2)$.
Let $\Pi(\theta)$ be the collection of $(P,Q)$ satisfying
Assumptions \ref{A1:strong:density}-\ref{A1:relative:signal},
where $\theta=\{\alpha, \beta_P, \beta_Q, \gamma, c_\lambda, r_\lambda, C_\beta, C_\alpha, C_\gamma\}$ is the collection of constants in the statements of the above assumptions. Based on the above assumptions and notation, the following theorem provides upper bounds for the excess risk of $\widehat{f}_{NN}$ and $\widehat{f}_{pa}$ in the special case $n_Q=0$.
\begin{theorem}\label{theorem:upper:bound:regret}
The following statements hold when $n_Q=0$:
\begin{enumerate}[label=(\alph*)]
\item \label{theorem:upper:bound:regret:smooth:source} If $\beta_P>\gamma \beta_Q$ and $k_P\asymp n_P^{\frac{2\beta_P}{2\beta_P+d}}$, then $\sup_{(P,Q)\in\Pi(\theta)}\mcR_Q(\widehat{f}_{NN})\lesssim n_P^{-\frac{(1+\alpha)\beta_P}{\gamma(2\beta_P+d)}}$;
\item  \label{theorem:upper:bound:regret:smooth:target} If $\beta_P\leq \gamma \beta_Q$ and $k_P\asymp n_P^{\frac{2\gamma\beta_Q}{2\gamma\beta_Q+d}}$, then $\sup_{(P,Q)\in\Pi(\theta)}\mcR_Q(\widehat{f}_{NN})\lesssim n_P^{-\frac{(1+\alpha)\beta_Q}{2\gamma\beta_Q+d}}$.
\end{enumerate}
Moreover, the following holds for the  classifier proposed in Algorithm \ref{alg:ag1}: 
\begin{enumerate}[label=(\alph*)]
\item If $\beta_P>\gamma \beta_Q$, then $\sup_{(P,Q)\in\Pi(\theta)}\mcR_Q(\widehat{f}_{\textrm{pa}})\lesssim \left(\frac{n_P}{\log^2(n_P)}\right)^{-\frac{(1+\alpha)\beta_P}{\gamma(2\beta_P+d)}}$;
\item If $\beta_P\leq \gamma \beta_Q$, then $\sup_{(P,Q)\in\Pi(\theta)}\mcR_Q(\widehat{f}_{\textrm{pa}})\lesssim \left(\frac{n_P}{\log^2(n_P)}\right)^{-\frac{(1+\alpha)\beta_Q}{2\gamma\beta_Q+d}}$.
\end{enumerate}
\end{theorem}
Theorem \ref{theorem:upper:bound:regret} provides upper bounds for the excess risk of $\widehat{f}_{NN}$ and $\widehat{f}_{pa}$ over $(P,Q)\in\Pi(\theta)$ in two smoothness scenarios: $\beta_P>\gamma\beta_Q$ and $\beta_P\le\gamma\beta_Q$. 
Up to logarithmic sacrifice, $\widehat{f}_{pa}$ performs equally well as $\widehat{f}_{NN}$.
When $\gamma$ is small, all upper bounds become smaller, indicating that more information has been transferred from $P$ to $Q$ to boost the classification performance.

Under certain circumstances, the excess risk has a very fast convergence rate. For instance, the excess risk is faster than $n_P^{-1/2}$ if $\beta_P>\gamma\beta_Q$
and $2(1+\alpha-\gamma)\beta_P\geq \gamma d$, or $\beta_P\leq \gamma \beta_Q$
and $2(1+\alpha-\gamma)\beta_Q\geq d$; it is faster than $n_P^{-1}$
if $\beta_P>\gamma\beta_Q$ and $(1+\alpha-2\gamma)\beta_P\geq \gamma d$,
or $\beta_P\le\gamma\beta_Q$ and $2(1+\alpha-\gamma)\beta_Q\geq d$.
It is easy to see that these results degenerate to \cite{audibert2007fast} 
in the conventional setting $P=Q$ and $\gamma=1$.
Our findings are summarized in Figure \ref{fig:phase_transition}(b). 

The following theorem provides the minimax lower bounds for the excess risk.
\begin{theorem}\label{theorem:lower:bound:regret}
If $n_Q=0$, then the following statements hold:
\begin{enumerate}[label=(\alph*)]
\item If $\beta_P>\gamma \beta_Q$ and $\alpha \beta_P\leq \gamma d$, then $\inf_{\widehat{f}}\sup_{(P, Q)\in \Pi(\theta)}\mcR_Q(\widehat{f})\geq cn_P^{-\frac{(1+\alpha)\beta_P}{\gamma(2\beta_P+d)}}$;
\item If $\beta_P\leq \gamma \beta_Q$ and $\alpha \beta_Q\leq d$, then $\inf_{\widehat{f}}\sup_{(P, Q)\in \Pi(\theta)}\mcR_Q(\widehat{f})\geq cn_P^{-\frac{(1+\alpha)\beta_Q}{2\gamma\beta_Q+d}}$,
\end{enumerate}
where $c$ is a positive constant relying on $\theta$, and the infimum is taken over classifiers constructed on the $P$-data.
\end{theorem}
We emphasize that the conditions $\alpha \beta_P\leq \gamma d$
and $\alpha\beta_Q\leq d$ are necessary to obtain minimax lower bounds. In fact, 
in the special case $P=Q$, we have $\beta_P=\beta_Q$ and $\gamma=1$, so both conditions reduce to $\alpha \beta_Q\leq d$, which was used by \cite{audibert2007fast} to establish the minimax lower bounds for the excess risk. Without assuming $\alpha \beta_Q \leq d$, the minimax lower bound remains unknown (see \citealp{pmlr-v65-locatelliandrea17a}).

Combining Theorem \ref{theorem:upper:bound:regret} and Theorem \ref{theorem:lower:bound:regret}, we immediately have the following conclusion:
\begin{enumerate}[label=(\alph*)]
\item If $\beta_P>\gamma \beta_Q$ and $\alpha \beta_P\leq \gamma d$, then $\inf_{\widehat{f}}\sup_{(P, Q)\in \Pi(\theta)}\mcR_Q(\widehat{f})\asymp n_P^{-\frac{(1+\alpha)\beta_P}{\gamma(2\beta_P+d)}}$;
\item If $\beta_P\leq \gamma \beta_Q$ and $\alpha \beta_Q\leq d$, then $\inf_{\widehat{f}}\sup_{(P, Q)\in \Pi(\theta)}\mcR_Q(\widehat{f})\asymp n_P^{-\frac{(1+\alpha)\beta_Q}{2\gamma\beta_Q+d}}$.
\end{enumerate}
Consequently, $\widehat{f}_{NN}$ and $\widehat{f}_{pa}$ both achieve the minimax optimal convergence rate. The optimal rate 
does not change when $\beta_P\in[0,\gamma\beta_Q]$, 
while it tends to zero faster when $\beta_Q\in(\gamma\beta_Q,1]$.
The above conclusions are summarized in Figure \ref{fig:phase_transition}(a) in which a phase transition phenomenon is observed.

\subsection{Extensions to General $n_Q$}\label{secton:two:source:model}
We extend the results in Section \ref{section:one:source} to general $n_Q$. We need the following assumption, a stronger version of Assumption \ref{A1:relative:signal}.
\begin{Assumption}\label{A1:relative:signal:multipe:source}(Mutual Relative Signal Exponent Condition) For all $\bfx \in \Omega$, it holds that
\begin{enumerate}[label=(\alph*)]
\item \label{A1:relative:signal:multipe:source:a} $(\eta_P(\bfx)-1/2)(\eta_Q(\bfx)-1/2)\geq 0$;
\item \label{A1:relative:signal:multipe:source:b}$ C_\gamma|\eta_Q(\bfx)-1/2|^\gamma \leq |\eta_P(\bfx)-1/2|\leq C_\gamma^{-1}|\eta_Q(\bfx)-1/2|^\gamma $ for some constants $\gamma, 0<C_\gamma<1$.
\end{enumerate}
\end{Assumption}
Assumption \ref{A1:relative:signal:multipe:source} requires an additional upper bound  $|\eta_P(\bfx)-1/2|\leq C_\gamma^{-1}|\eta_Q(\bfx)-1/2|^\gamma$ in comparison with Assumption \ref{A1:relative:signal}. After rewriting this bound as $|\eta_Q(\bfx)-1/2|\geq C_\gamma^{1/\gamma}|\eta_P(\bfx)-1/2|^{1/\gamma}$, we can see the additional condition essentially requires $1/\gamma$ being the relative signal exponent of $Q$ with respective to $P$, which assesses the  information that can be transferred from $Q$ to $P$ (see the discussion right after Theorem \ref{theorem:upper:bound:regret}). Therefore, it is reasonable to view $\gamma$ as the \textit{mutual relative signal exponent}.

Let $\Pi'(\theta)$ be the collection of $(P, Q)$ satisfying Assumptions \ref{A1:strong:density}, \ref{A1:holder:smooth}, \ref{A1:marginal:assumption} and \ref{A1:relative:signal:multipe:source}, where
$\theta=\{\alpha, \beta_P, \beta_Q, \gamma, c_\lambda, r_\lambda, C_\beta, C_\alpha, C_\gamma\}$ is the collection of constants involved in the corresponding assumptions.
Clearly, $\Pi'(\theta)$ is a subset of $\Pi(\theta)$.
 The following theorem extends the results in Theorem \ref{theorem:upper:bound:regret} to general $n_Q$.
\begin{theorem}\label{theorem:two:sample:KNN:regret}
Suppose that either $n_Q\to\infty$ or $n_P\to\infty$. Then the following statements hold:
\begin{enumerate}[label=(\alph*)]
\item If $\beta_P>\gamma \beta_Q$, $w_Q\asymp\delta$, $w_P\asymp\delta^\gamma, k_Q\asymp n_Q\delta^{\frac{\gamma d}{\beta_P}}$, and $k_P\asymp n_P\delta^{\frac{\gamma d}{\beta_P}}$, where $\delta=(n_P^{\frac{2\beta_P+\gamma d}{\gamma (2\beta_P+d)}}+n_Q)^{-\frac{\beta_P}{2\beta_P+\gamma d}}$, then
\begin{equation}\label{eq:theorem:two:sample:KNN:regret:source} 
\sup_{(P,Q)\in\Pi'(\theta)}\mcR_Q(\widehat{f}_{NN})\lesssim (n_P^{\frac{2\beta_P+\gamma d}{\gamma(2\beta_P+d)}}+n_Q)^{-\frac{\beta_P(1+\alpha)}{2\beta_P+\gamma d}};
\end{equation}
\item If $\beta_P\leq\gamma \beta_Q$, 
$w_Q\asymp\delta,w_P\asymp\delta^\gamma, k_Q\asymp n_Q\delta^{\frac{d}{\beta_Q}}$, and $k_P\asymp n_P\delta^{\frac{ d}{\beta_Q}}$,
where $\delta=(n_P^{\frac{2\beta_Q+d}{2\gamma\beta_Q+d}}+n_Q)^{-\frac{\beta_Q}{2\beta_Q+d}}$, then 
\begin{equation}\label{eq:theorem:two:sample:KNN:regret:target}
\sup_{(P,Q)\in\Pi(\theta)}\mcR_Q(\widehat{f}_{NN})\lesssim (n_P^{\frac{2\beta_Q+d}{2\gamma\beta_Q+d}}+n_Q)^{-\frac{\beta_Q(1+\alpha)}{2\beta_Q+d}}.
\end{equation}
\end{enumerate}
Moreover, the  adaptive classifier proposed in Algorithm \ref{alg:ag2} has the following properties: 
\begin{enumerate}[label=(\alph*)]
\item If $\beta_P>\gamma \beta_Q$, then
\begin{equation}
\sup_{(P,Q)\in\Pi'(\theta)}\mcR_Q(\widehat{f}_{pa})\lesssim \bigg[\bigg(\frac{n_P}{\log^2(n_P+n_Q)}\bigg)^{\frac{2\beta_P+\gamma d}{\gamma(2\beta_P+d)}}+\frac{n_P}{\log^2(n_P+n_Q)}\bigg]^{-\frac{\beta_P(1+\alpha)}{2\beta_P+\gamma d}};\nonumber
\end{equation}
\item If $\beta_P\leq\gamma \beta_Q$, then
\begin{equation}
\sup_{(P,Q)\in\Pi(\theta)}\mcR_Q(\widehat{f}_{pa})\lesssim \bigg[\bigg(\frac{n_P}{\log^2(n_P+n_Q)}\bigg)^{\frac{2\beta_Q+d}{2\gamma\beta_Q+d}}+\frac{n_Q}{\log^2(n_P+n_Q)}\bigg]^{-\frac{\beta_Q(1+\alpha)}{2\beta_Q+d}}.\nonumber
\end{equation}
\end{enumerate}
\end{theorem}
Theorem \ref{theorem:two:sample:KNN:regret} implies that
the excess risks of $\widehat{f}_{NN}$ and $\widehat{f}_{pa}$ 
have the same 
upper bounds under general $n_Q$.
The bounds under $\beta_P>\gamma\beta_Q$ are derived over $\Pi'(\theta)$,
which is a subset of  $\Pi(\theta)$,
so it is interesting to explore the upper bounds 
over $\Pi(\theta)$ as well. 
A reexmination of the proof reveals that, if $\beta_P> \gamma \beta_Q$, then
\begin{equation}
\sup_{(P,Q)\in\Pi(\theta)}\mcR_Q(\widehat{f}_{NN})\lesssim 
 (n_P^{\frac{(2\beta_Q+d)\beta_P}{\gamma(2\beta_P+d)\beta_Q}}+n_Q)^{-\frac{\beta_Q(1+\alpha)}{2\beta_Q+ d}}.
\label{eq:suboptimal:upper:bound}
\end{equation}
When $n_Q=0$, (\ref{eq:theorem:two:sample:KNN:regret:target}) and (\ref{eq:suboptimal:upper:bound}) degenerate to
Theorem \ref{theorem:upper:bound:regret} both being optimal thanks to Theorem \ref{theorem:lower:bound:regret}.
For general $n_Q$, (\ref{eq:theorem:two:sample:KNN:regret:target}) is optimal thanks to Theorem \ref{theorem:lowerbound:multiple:sources} \ref{eq:theorem:two:sample:KNN:regret:target:b}, whereas (\ref{eq:suboptimal:upper:bound}) is substantially slower than the lower bound stated in Theorem \ref{theorem:lowerbound:multiple:sources} \ref{eq:theorem:two:sample:KNN:regret:target:a}.

\begin{theorem}\label{theorem:lowerbound:multiple:sources}
The following statements hold:
\begin{enumerate}[label=(\alph*)]
\item\label{eq:theorem:two:sample:KNN:regret:target:a} If $\beta_P>\gamma \beta_Q$ and  $\alpha \beta_P\leq \gamma d$, then 
\begin{eqnarray}
\inf_{\widehat{f}}\sup_{(P, Q)\in \Pi(\theta)}\mcR_Q(\widehat{f})\geq \inf_{\widehat{f}}\sup_{(P, Q)\in \Pi'(\theta)}\mcR_Q(\widehat{f})\geq c(n_P^{\frac{2\beta_P+\gamma d}{\gamma(2\beta_P+d)}}+n_Q)^{-\frac{\beta_P(1+\alpha)}{2\beta_P+\gamma d}};\nonumber
\end{eqnarray}
\item\label{eq:theorem:two:sample:KNN:regret:target:b} If $\beta_P\leq \gamma \beta_Q$ and $\alpha \beta_Q\leq d$, then 
\begin{eqnarray}
\inf_{\widehat{f}}\sup_{(P, Q)\in \Pi(\theta)}\mcR_Q(\widehat{f})\geq \inf_{\widehat{f}}\sup_{(P, Q)\in \Pi'(\theta)}\mcR_Q(\widehat{f})\geq c (n_P^{\frac{2\beta_Q+d}{2\gamma\beta_Q+d}}+n_Q)^{-\frac{\beta_Q(1+\alpha)}{2\beta_Q+d}},\nonumber
\end{eqnarray}
\end{enumerate}
where $c$ is a constant depending on $\theta$, and the infimum is taken over the classifiers constructed based on the entire data $\mcZ$ described in (\ref{eqn:entire:data}).
\end{theorem}

Theorem \ref{theorem:lowerbound:multiple:sources} provides lower bounds
for the minimax excess risk under $\Pi(\theta)$ and $\Pi'(\theta)$. 
Combining Theorem \ref{theorem:two:sample:KNN:regret} and Theorem
\ref{theorem:lowerbound:multiple:sources}, we get that
\begin{enumerate}[label=(\alph*)]
\item \textrm{if $\beta_P>\gamma \beta_Q$ and $\alpha \beta_P\leq \gamma d$},\textrm{ then $\inf_{\widehat{f}}\sup_{(P, Q)\in \Pi'(\theta)}\mcR_Q(\widehat{f})\asymp 
(n_P^{\frac{2\beta_P+\gamma d}{\gamma(2\beta_P+d)}}+n_Q)^{-\frac{\beta_P(1+\alpha)}{2\beta_P+\gamma d}}$;}
\item \textrm{if $\beta_P\leq \gamma \beta_Q$ and $\alpha \beta_Q\leq d$},\textrm{ then
$\inf_{\widehat{f}}\sup_{(P, Q)\in \Pi'(\theta)}\mcR_Q(\widehat{f})\asymp 
(n_P^{\frac{2\beta_Q+d}{2\gamma\beta_Q+d}}+n_Q)^{-\frac{\beta_Q(1+\alpha)}{2\beta_Q+d}}$.}
\end{enumerate}
In view of Theorems \ref{theorem:upper:bound:regret}-\ref{theorem:lowerbound:multiple:sources} and (\ref{eq:suboptimal:upper:bound}), we summarize the (sub)optimality of convergence rate for the minimax excess risk in the following Table \ref{table:con:diagram}. 
\begin{table}[H]
\caption{\it (Sub)optimality of convergence rate for the minimax excess risk in different regimes.}
\begin{tabular}{cccc}
 &  & $\beta_P \leq \gamma \beta_Q, \alpha \beta_Q\leq d$ & $\beta_P>\gamma\beta_Q, \alpha\beta_P\leq \gamma d$ \\ \hline
\multirow{2}{*}{$n_Q=0$} & rate for $\inf_{\widehat{f}}\sup_{(P, Q)\in \Pi'(\theta)}\mcR_Q(\widehat{f})$ & Optimal & Optimal \\ 
& rate for $\inf_{\widehat{f}}\sup_{(P, Q)\in \Pi(\theta)}\mcR_Q(\widehat{f})$ & Optimal & Optimal \\ \hline
\multirow{2}{*}{$n_Q>0$} & rate for $\inf_{\widehat{f}}\sup_{(P, Q)\in \Pi'(\theta)}\mcR_Q(\widehat{f})$ & Optimal & Optimal \\ 
& rate for $\inf_{\widehat{f}}\sup_{(P, Q)\in \Pi(\theta)}\mcR_Q(\widehat{f})$ & Optimal & Suboptimal \\ \hline
\end{tabular}
\label{table:con:diagram}
\end{table}

%\begin{eqnarray}
%	 \delta=(n_P^{\frac{\beta_P(2\beta_Q+d)}{\gamma(2\beta_P+d)\beta_Q}}+n_Q)^{-\frac{\beta_Q}{2\beta_Q+d}},  w_Q=\delta,w_P=\delta^\gamma, k_Q=\floor{n_Q\delta^{\frac{d}{\beta_Q}}}, k_P=\floor{n_P\delta^{\frac{\gamma d}{\beta_P}}}.\nonumber
%\end{eqnarray}
\subsection{Extensions to Multiple Sources}\label{section:multiple:source:model}
In this section, we extend the previous results to the scenario where the data come from multiple sources. Multi-source scenario is common in big data research (\citealp{zhang2015divide, shang2017computational}).  Suppose that, for $1\le j\le m$, the observations $(\bfX_1^{P_j}, Y_1^{P_j}), \ldots, (\bfX_{n_j}^{P_j}, Y_{n_j}^{P_j})$ are generated from a source distribution $P_j$. Without loss of generality, assume $n_1\geq n_2\geq \ldots \geq n_m$. 
Let $Q$ be the target distribution. For simplicity, assume $n_Q=0$,
though the results are extendable to general $n_Q$. Similar to Section \ref{section:one:source}, define the $k$NN classifier as follows:
\[
\widehat{f}_{NN}(\bfx)=\mbI(\widehat{\eta}_{NN}(\bfx)\geq 1/2),
\]
where $\widehat{\eta}_{NN}(\bfx)=\sum_{j=1}^mw_jk_j\widehat{\eta}_{k_j}(\bfx)/
\sum_{j=1}^mw_jk_j$,
$\widehat{\eta}_{k_j}(\bfx)$ is the $k$NN estimator of $P_j(Y^{P_j}=1|\bfX^{P_j}=\bfx)$ based on the $k_j$ nearest covariates in $P_j$-data, and $w_j>0$ is the corresponding weight. We also propose an adaptive classifier $\widehat{f}_{pa}$ in Algorithm \ref{alg:ag3} in which the parameters $w_j$'s and $k_j$'s are data-driven. 
Note that Algorithm \ref{alg:ag3} selects the tuple $(k_1,\ldots,k_m)$ over
$\{(k_1,\ldots,k_m): k_j=\floor{k_1n_j/n_1}, 2\le j\le m, 1\le k_1\le n_1\}$,
which requires $\max\{n_1,\ldots,n_m\} (\textrm{which is } n_1)$ attempts. 
In contrast, \cite{cai2019transfer} requires $n_1+\cdots+n_m$ attempts to select the
tuple, hence, the ratio of attempts for 
\cite{cai2019transfer} and Algorithm \ref{alg:ag3} is $\frac{n_1+\cdots+n_m}{\max\{n_1,\ldots,n_m\}}$
which is nearly $m$ if $n_1=\cdots=n_m$.

%The proposed algorithms scan the nearest covariates in local data, and so, are computationally friendly. Similar to Section \ref{section:one:source}, we can show that ratio of the computational costs between the adaptive algorithm in \cite{cai2019transfer} and Algorithm \ref{alg:ag3} is $\sum_{j=1}^mn_j/\max_{1\leq j\leq m}n_j$, which is greater than one and even achieves $m$ when all $n_j$'s are equal.

\begin{algorithm}[htp!]
\SetAlgoLined
\KwIn{data $( \bfX_1^{P_j}, Y_1^{P_j}),\ldots, (\bfX_{n_{j}}^{P_j}, Y_{n_j}^{P_j})$ for $j=1,\ldots, m$ and new features $\bfx$;}
% \textbf{Initialization:} $k=1$\;
\textbf{Initiation: } set $k_1=0$;

\textbf{while} $k_1< n_1$ \textbf{do}

\hspace*{0.5cm}  update $k_1:=k_1+1$ and $k_j:=\floor{k_1n_j/n_1}$ for $j=2,\ldots, m$;

\hspace*{0.5cm} calculate $\widehat{\eta}_{k_j}(\bfx)=\frac{1}{k_j}\sum_{i=1}^{k_j} Y_{(i)}^{P_j}(\bfx)$ (set $\widehat{\eta}_{k_j}=1/2$ if $k_j=0$) for $j=1,\ldots, m$;

\hspace*{0.5cm} calculate $r_{k_1}^+=\sqrt{\sum_{j=1}^m\mbI(\widehat{\eta}_{k_j}(\bfx)\geq 1/2)k_j(\widehat{\eta}_{k_j}(\bfx)-1/2)^2}$

\hspace*{0.5cm} calculate $r_{k_1}^-=\sqrt{\sum_{j=1}^m\mbI(\widehat{\eta}_{k_j}(\bfx)<1/2)k_j(\widehat{\eta}_{k_j}(\bfx)-1/2)^2}$

\hspace*{0.5cm} calculate  $r_{k_1}=\max(r_{k_1}^+, r_{k_1}^-)$;

\hspace*{0.5cm} \textbf{if } $r_{k_1}>\sqrt{[d+\log(\sum_{j=1}^sn_j)]\log(\sum_{j=1}^sn_j)}$   \textbf{or} $k_1=n_1$ \textbf{ then} 

\hspace*{1cm}  set $\widehat{k}_j=k_j$ for $j=1,\ldots, m$;

\hspace*{1cm}  exit  loop;

\hspace*{0.5cm} \textbf{end if}

\textbf{end while}

calculate $\widehat{\eta}_{\widehat{k}_j}(\bfx)$ (set $\widehat{\eta}_{\widehat{k}_j}(\bfx)=1/2$ if $\widehat{k}_j=0$) for $j=1,\ldots, m$;

%calculate $\widehat{\eta}_{pa}(\bfx)=\sum_{j=1}^m\widehat{\lambda}_j \widehat{\eta}_{\widehat{k}_j}(\bfx)$;

\KwOut{classifier $\widehat{f}_{pa}(\bfx)=\mbI(\sum_{j=1}^mk_j(\widehat{\eta}_{\widehat{k}_j}(\bfx)-1/2)\geq 0)$.}
 \caption{Multiple-Sample Pointwise Adaptive KNN}
 \label{alg:ag3}
\end{algorithm}

We extend the theoretical results in Sections \ref{section:one:source} and \ref{secton:two:source:model} to multi-source scenario. For that,
let $\theta_m=\{\alpha, \beta_1, \ldots, \beta_m, \beta_Q, \gamma_1,\ldots, \gamma_m, c_\lambda, r_\lambda, C_\beta, C_\alpha, C_\gamma\}$
and
$\Pi_m'(\theta_m)$ be the collection of tuples $(P_1, \ldots, P_m, Q)$ such that $(P_j, Q)\in \Pi'(\alpha, \beta_j, \beta_Q, \gamma_j, c_\lambda, r_\lambda, C_\beta, C_\alpha, C_\gamma)$ for $j=1,\ldots, m$.
\begin{theorem}\label{theorem:multiple:sample:KNN:regret}
Let $\beta^*=\min\{\frac{\beta_1}{\gamma_1}, \ldots,\frac{\beta_m}{\gamma_m}, \beta_Q\}$, $\delta=(\sum_{s=1}^mn_s^{\frac{2\beta^*+d}{2\gamma_s\beta^*+d}})^{-\frac{\beta^*}{2\beta^*+d}}$. If $w_j\asymp \delta^{\gamma_j}$ and $k_j\asymp n_j\delta^{\frac{d}{\beta^*}}$ for all $j=1,\ldots, m$, then the following holds:
\begin{eqnarray}\label{eqn:multi-source}
&&{\sup_{(P_1,\ldots,P_m,Q)\in\Pi'_m(\theta_m)}}\mcR_Q(\widehat{f}_{NN})\lesssim \bigg(\sum_{s=1}^mn_s^{\frac{2\beta^*+d}{2\gamma_s\beta^*+d}}\bigg)^{-\frac{\beta^*(1+\alpha)}{2\beta^*+d}},\\ &&{\sup_{(P_1,\ldots,P_m,Q)\in\Pi'_m(\theta_m)}}\mcR_Q(\widehat{f}_{pa})\lesssim \bigg[\sum_{s=1}^m\bigg(\frac{n_s}{\log^2(\sum_{j=1}^mn_j)}\bigg)^{\frac{2\beta^*+d}{2\gamma_s\beta^*+d}}\bigg]^{-\frac{\beta^*(1+\alpha)}{2\beta^*+d}}.\nonumber
\end{eqnarray}
Furthermore, there exists a constant $c<0$ depending on $\theta_m$ such that
\begin{eqnarray}
\inf_{\widehat{f}}\sup_{(P_1, \ldots, P_m, {Q})\in \Pi'_m(\theta_m)}\mcR_Q(\widehat{f})\geq c  \bigg(\sum_{s=1}^mn_s^{\frac{2\beta^*+d}{2\gamma_s\beta^*+d}}\bigg)^{-\frac{\beta^*(1+\alpha)}{2\beta^*+d}},\nonumber
\end{eqnarray}
where the infimum is taken over the classifiers based on the entire data
$(\bfX_i^{P_j}, Y_i^{P_j})$, $j=1,\ldots,m$, $i=1,\ldots,n_j$.
\end{theorem}
Theorem \ref{theorem:multiple:sample:KNN:regret} derives an exact order for the minimax excess risk in multi-source scenario:
\[
\inf_{\widehat{f}}\sup_{(P_1, \ldots, P_m, {Q})\in \Pi'_m(\theta_m)}\mcR_Q(\widehat{f})\asymp\left(\sum_{s=1}^mn_s^{\frac{2\beta^*+d}{2\gamma_s\beta^*+d}}\right)^{-\frac{\beta^*(1+\alpha)}{2\beta^*+d}}.
\]
In the special case with $m=2$, $P_1=P$ and $P_2=Q$,   
the RHS of (\ref{eqn:multi-source})
becomes the upper bounds in Theorem \ref{theorem:two:sample:KNN:regret}. 
Moreover, both $\widehat{f}_{NN}$ and $\widehat{f}_{pa}$ are proven minimax optimal. The proof of Theorem \ref{theorem:multiple:sample:KNN:regret} is similar to Theorems \ref{theorem:two:sample:KNN:regret} and \ref{theorem:lowerbound:multiple:sources}, hence, is omitted.

\section{Monte Carlo Experiments}\label{section:simulation}
In this section, we investigate the finite-sample performance of the proposed algorithms through Monte Carlo experiments. 
We chose $P, Q$, the source and target distributions, 
to be uniform on $[0, 1]^2$. 
The data generating process (DGP) proceeds by first generating features $\bfX\sim Q$ and then generating label $Y\sim\eta_Q$. Various choices of $\eta_P,\eta_Q$ are summarized below.

\begin{itemize}[wide, labelwidth=!, labelindent=0pt]
\item[DGP 1:]  For $\bfx\in [0, 1]^2$, $\kappa\in [0, 1]$, $\gamma\in (0, \infty)$,
\begin{eqnarray*}
\eta_Q(\bfx)=\begin{cases}
\kappa\bigg(\frac{\|\bfx\|}{\sqrt{2}}-\frac{1}{2}\bigg)+\frac{1}{2}, &\textrm{if the ten-billionth value of $\|\bfx\|/\sqrt{2}$ is even};\\  
\kappa^\gamma\bigg(\frac{\|\bfx\|}{\sqrt{2}}-\frac{1}{2}\bigg)+\frac{1}{2}, &\textrm{if the ten-billionth value of $\|\bfx\|/\sqrt{2}$ is odd},\\  
\end{cases}
\end{eqnarray*}
and
\begin{eqnarray*}
\eta_P(\bfx)=sign\bigg(\frac{\|\bfx\|}{\sqrt{2}}-\frac{1}{2}\bigg)\kappa^\gamma\bigg|\frac{\|\bfx\|}{\sqrt{2}}-\frac{1}{2}\bigg|^\gamma +\frac{1}{2}.
\end{eqnarray*}

\item[DGP 2:]  For $\bfx \in [0, 1]^d$, $\kappa\in [0, 1]$, $\gamma\in (0, \infty)$,
\begin{eqnarray*}
\eta_Q(\bfx)=\kappa\bigg(\frac{\|\bfx\|}{\sqrt{2}}-\frac{1}{2}\bigg)+\frac{1}{2},
\end{eqnarray*}
and
\begin{eqnarray*}
\eta_P(\bfx)=\begin{cases}
sign\bigg(\frac{\|\bfx\|}{\sqrt{2}}-\frac{1}{2}\bigg) \kappa^\gamma \bigg|\frac{\|\bfx\|}{\sqrt{2}}-\frac{1}{2}\bigg|^\gamma+\frac{1}{2}, &\textrm{if the ten-billionth value of $\|\bfx\|/\sqrt{2}$ is even};\\  
sign\bigg(\frac{\|\bfx\|}{\sqrt{2}}-\frac{1}{2}\bigg)(1.2\kappa)^\gamma \bigg|\frac{\|\bfx\|}{\sqrt{2}}-\frac{1}{2}\bigg|^\gamma+\frac{1}{2}, &\textrm{if the ten-billionth value of $\|\bfx\|/\sqrt{2}$ is odd}.\\  
\end{cases}
\end{eqnarray*}
\end{itemize}
We will comment both DGPs satisfy Assumptions \ref{A1:strong:density},
\ref{A1:holder:smooth}, \ref{A1:marginal:assumption}, \ref{A1:relative:signal:multipe:source} are all satisfied.
Since $\bfX$ is uniformly distributed on $[0, 1]^2$ under both $P$ and $Q$, Assumption \ref{A1:strong:density} holds. Obviously, we have $(\beta_Q, \beta_P)=(0, \gamma)$ and $(\beta_Q, \beta_P)=(1, 0)$ which implies Assumption \ref{A1:holder:smooth}. We can verify Assumption \ref{A1:marginal:assumption}  with $\alpha=1$ in both DPGs, and show that  Assumption \ref{A1:relative:signal:multipe:source} is satisfied for both DGPs with 
\[
\kappa^{\gamma(1-\gamma)}|\eta_Q(\bfx)-1/2|^\gamma\leq |\eta_P(\bfx)-1/2|\leq |\eta_Q(\bfx)-1/2|^\gamma
\]
and
\[
|\eta_Q(\bfx)-1/2|^\gamma\leq |\eta_P(\bfx)-1/2|\leq 1.2^\gamma|\eta_Q(\bfx)-1/2|^\gamma.
\]
Heuristically, with a larger $\kappa$, both $|\eta_Q(\bfx)-1/2|$ and $|\eta_P(\bfx)-1/2|$ become stronger, which makes the classification problem easier. In both DGPs, we chose $\gamma=0.6$, $\kappa=0.1,0.2,\ldots, 0.9$ and $n_Q=2000$, $n_P=200, 500, 1000, 2000, 3500, 5000$.  We considered three competitors: ($k$NNCW) the adaptive algorithm proposed in \cite{cai2019transfer}; ($k$NNQ) naive $k$NN on $Q$-data only; ($k$NNALL) naive $k$NN on the entire data $\mcZ$ (recall that $\mcZ$ is the collection of both $P$-data and $Q$-data). 
To approximate the classification accuracy, we generated new features $\bfx^{\textrm{new}}$ and calculated their predicted label $y^*=\mbI(\eta_Q(\bfx^{\textrm{new}})\geq 1/2)$ based on Bayes classifier, and the classification accuracy of the proposed algorithms is approximated by the percentage of producing the same prediction as $y^*$ over $1000$ replicated trials. Moreover, we compare the runtime of Algorithm \ref{alg:ag2} and $k$NNCW in different settings.

\begin{figure}[h]
\centering

\includegraphics[width=2 in]{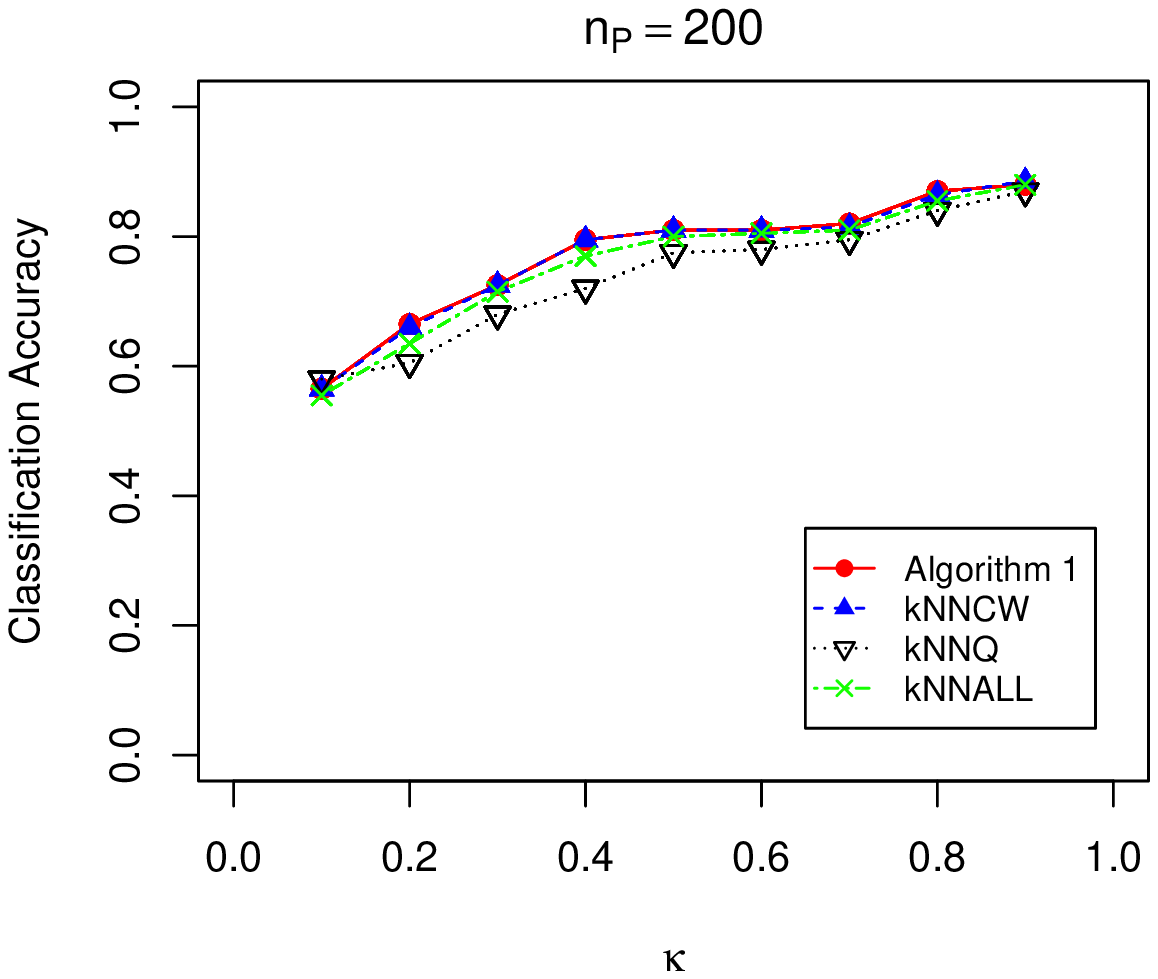}
\includegraphics[width=2 in]{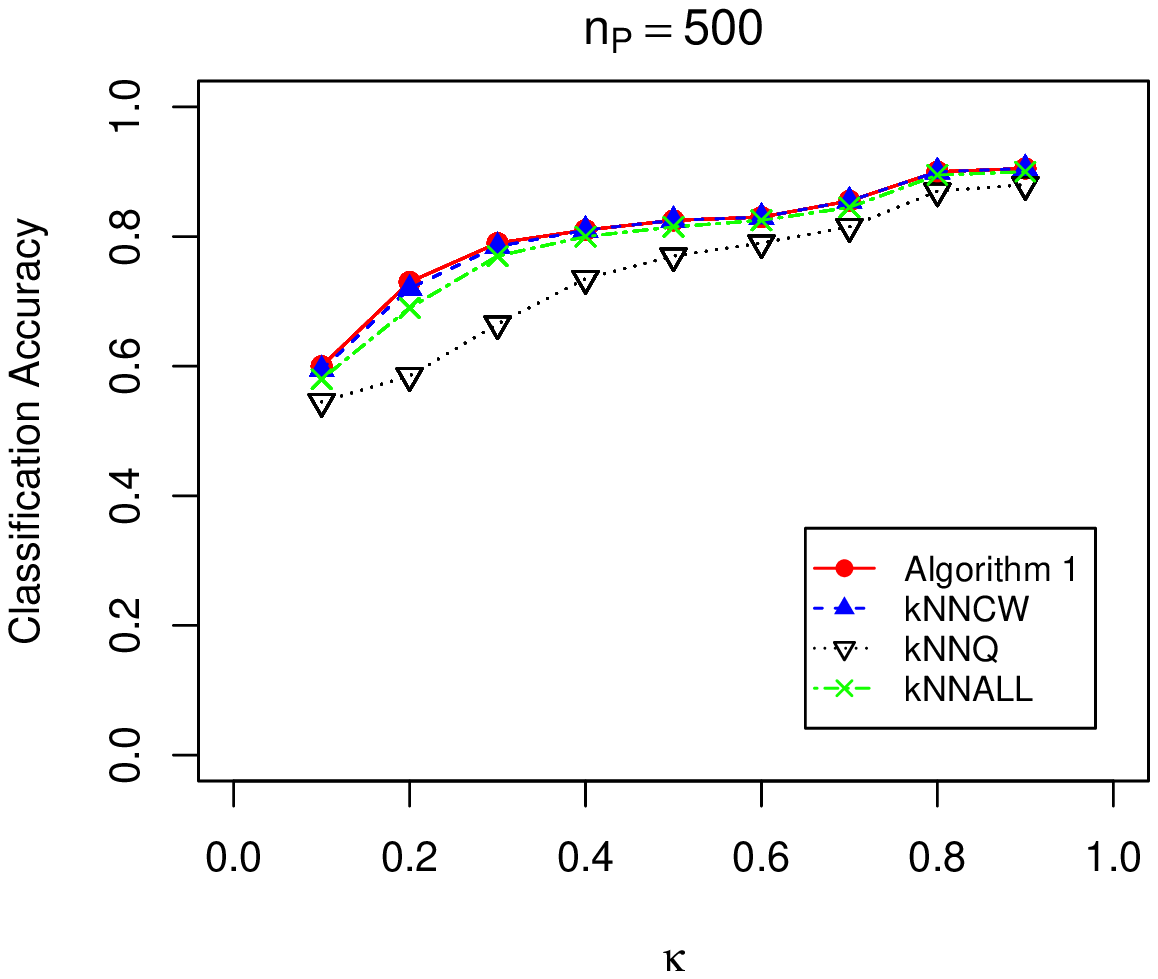}
\includegraphics[width=2 in]{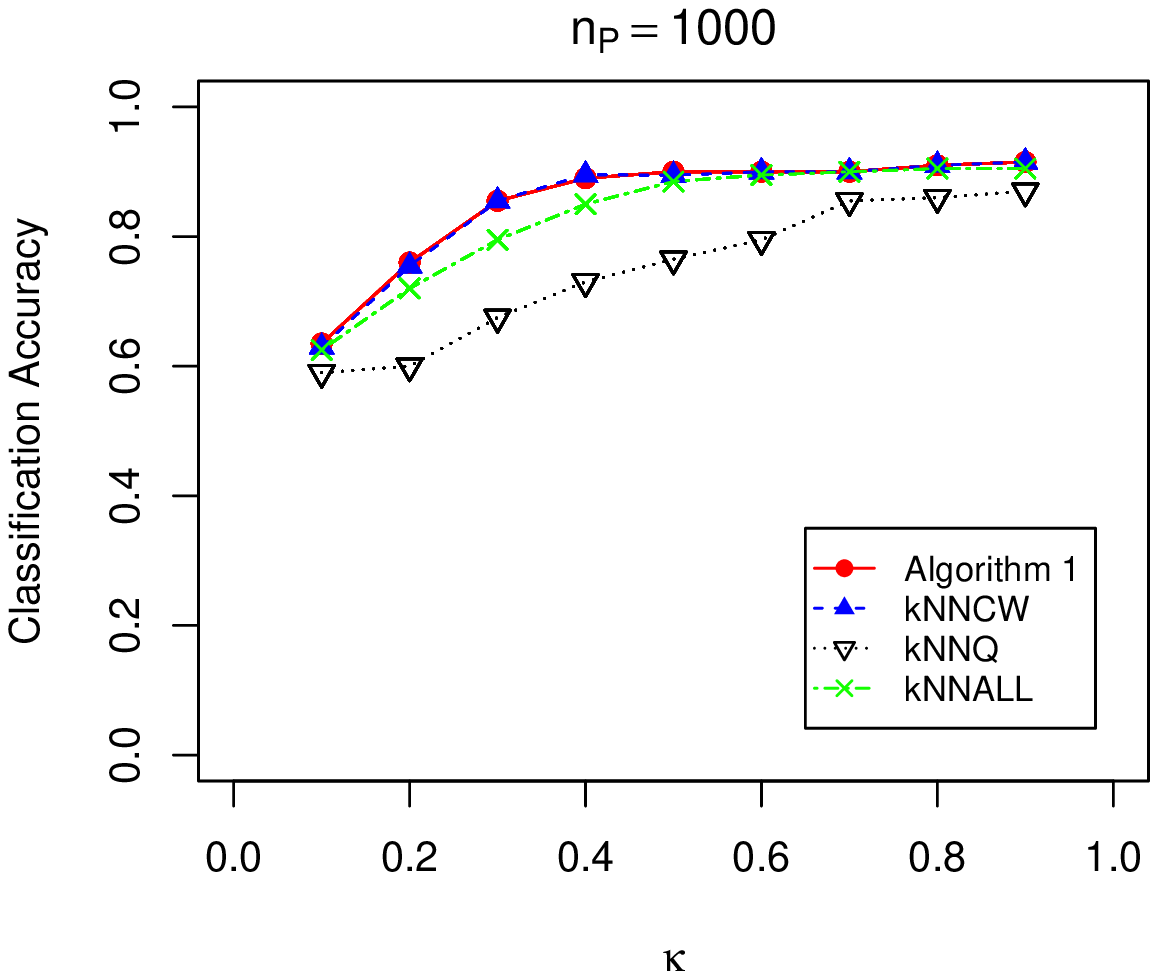}

\includegraphics[width=2 in]{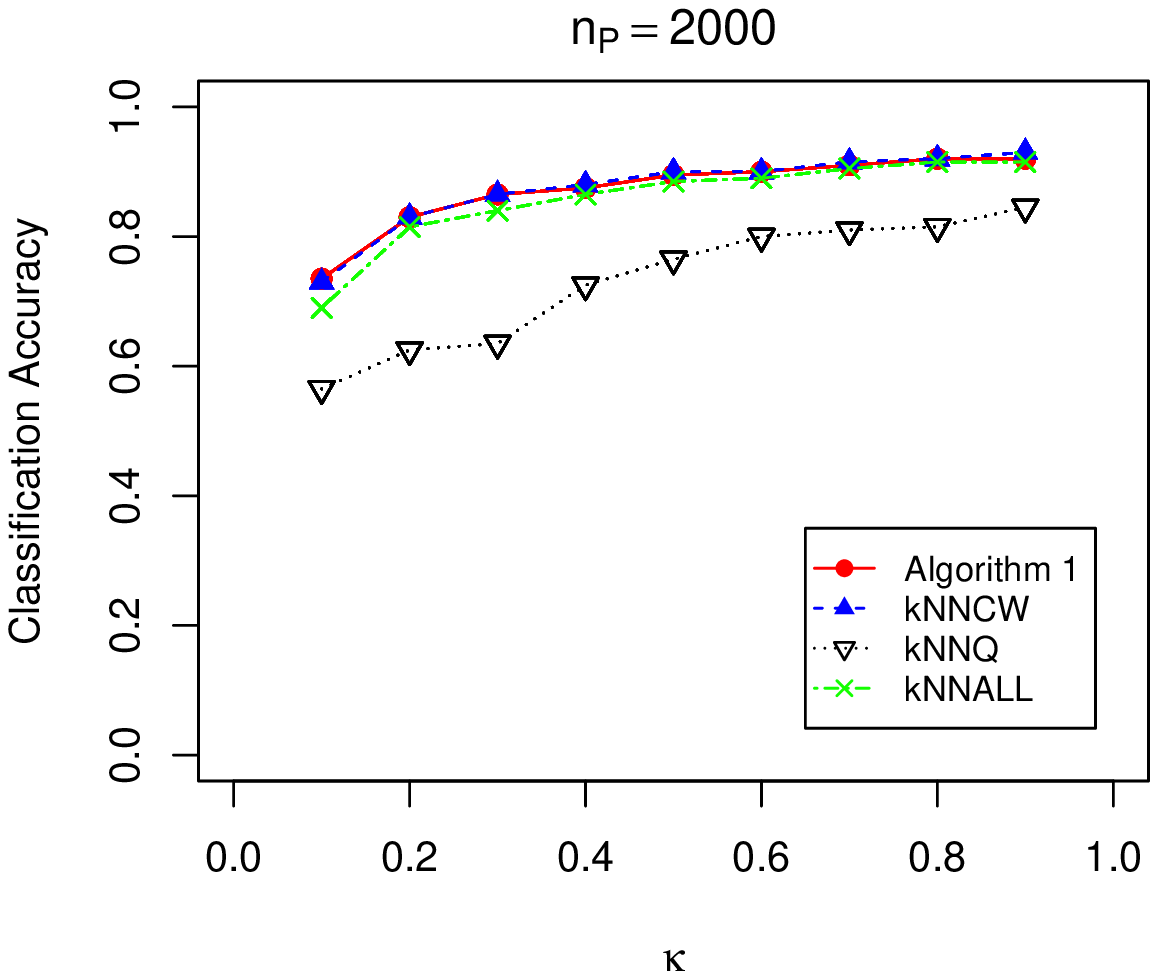}
\includegraphics[width=2 in]{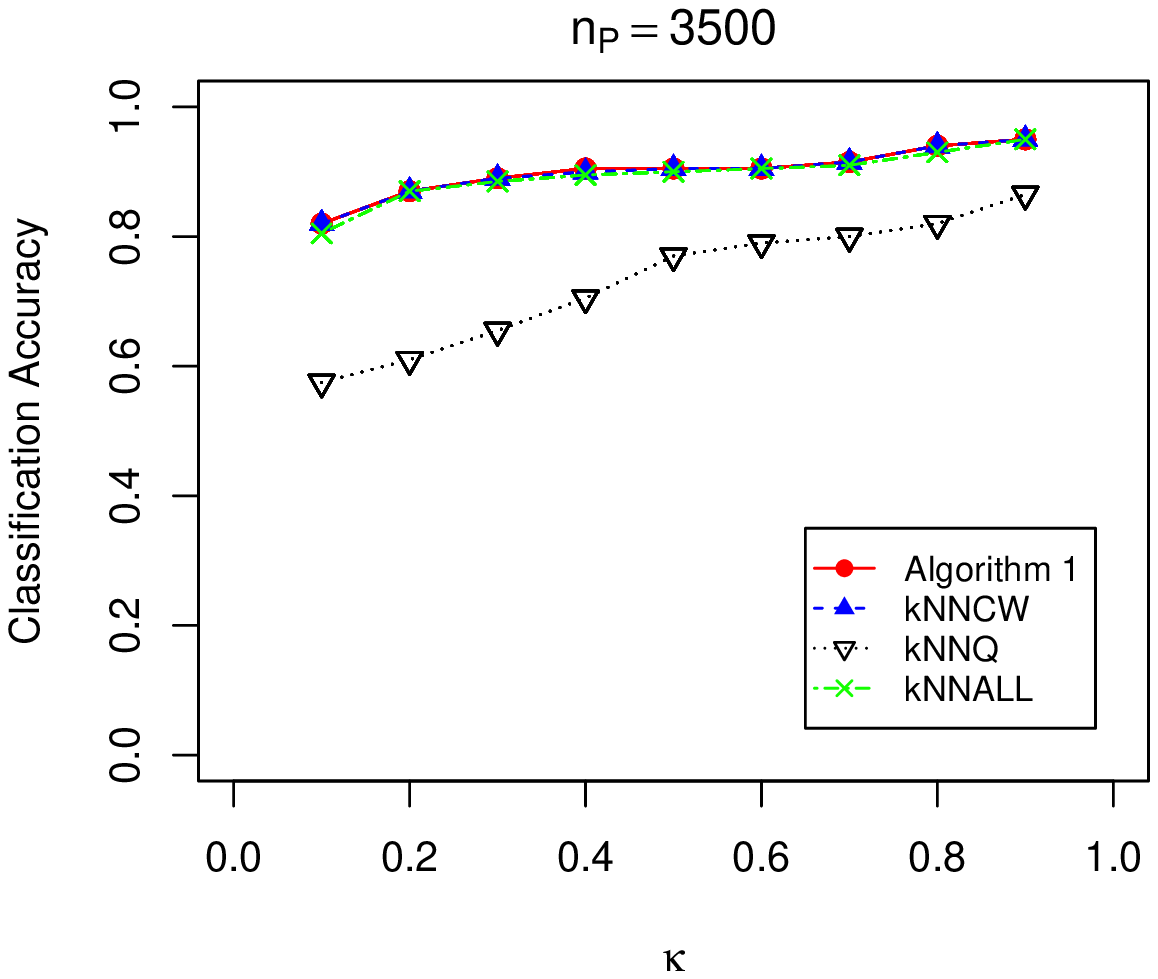}
\includegraphics[width=2 in]{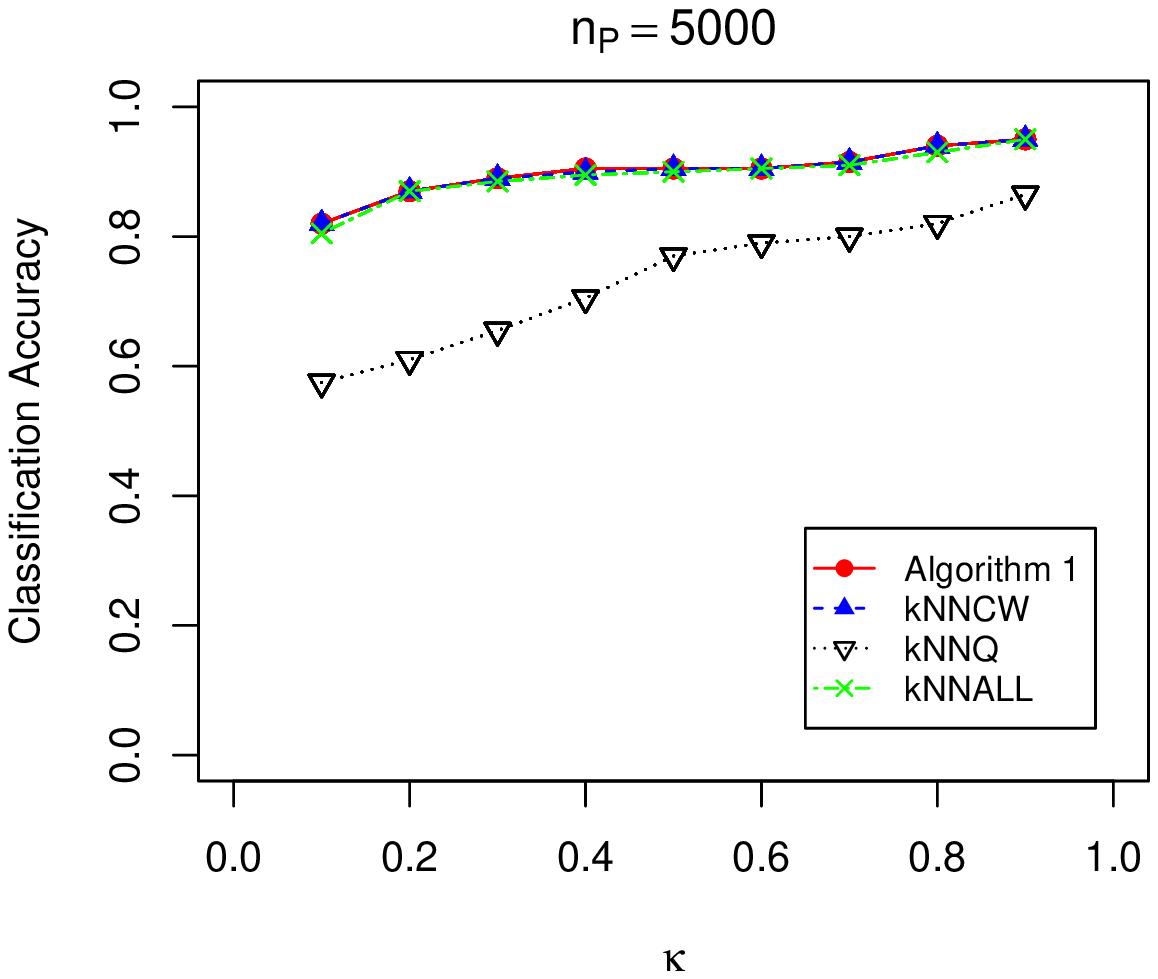}
\caption{\it DGP 1 (Smooth Source): Classification accuracy under different combinations of $(n_P, \kappa)$.}
\label{figure:classification:DGP1}
\end{figure}

\begin{figure}[h]
\centering

\includegraphics[width=2 in]{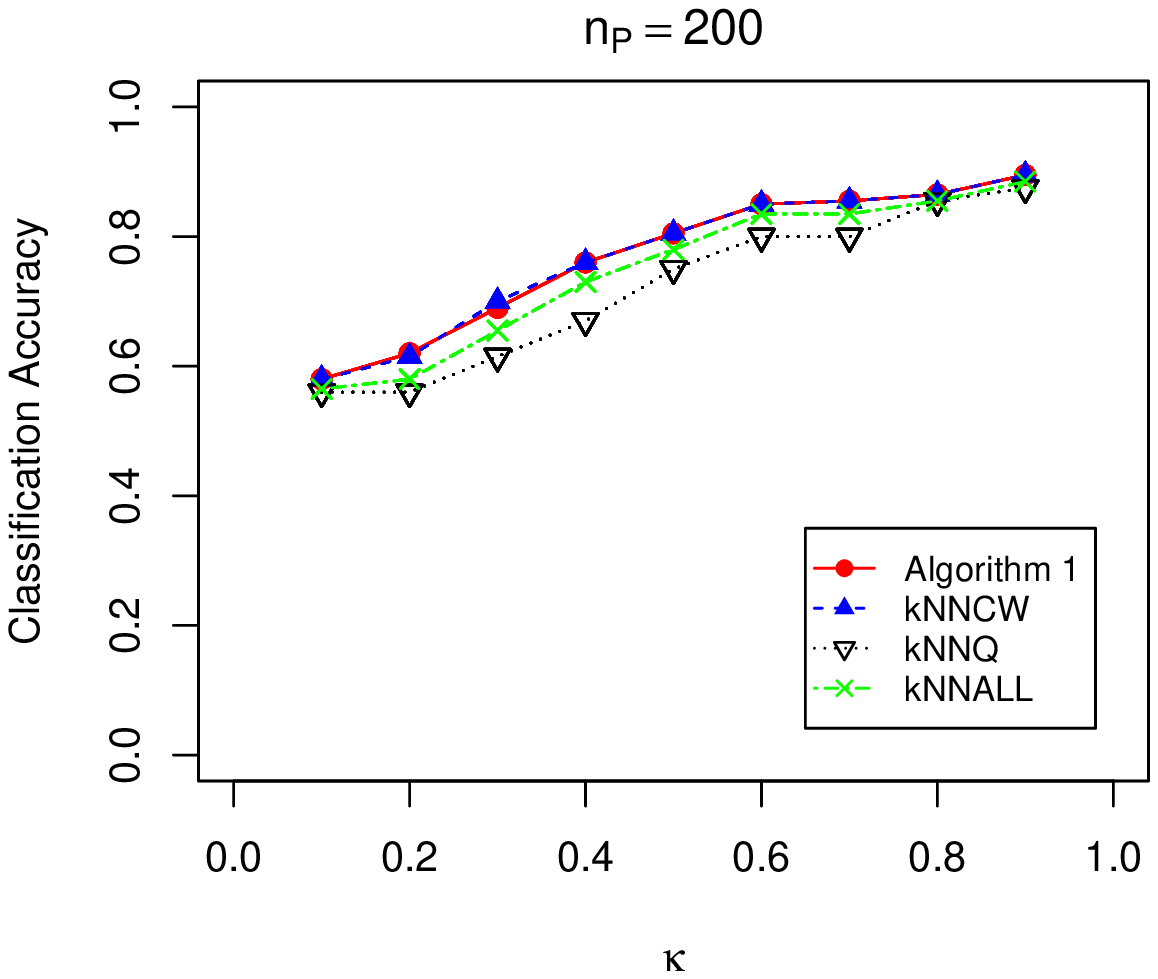}
\includegraphics[width=2 in]{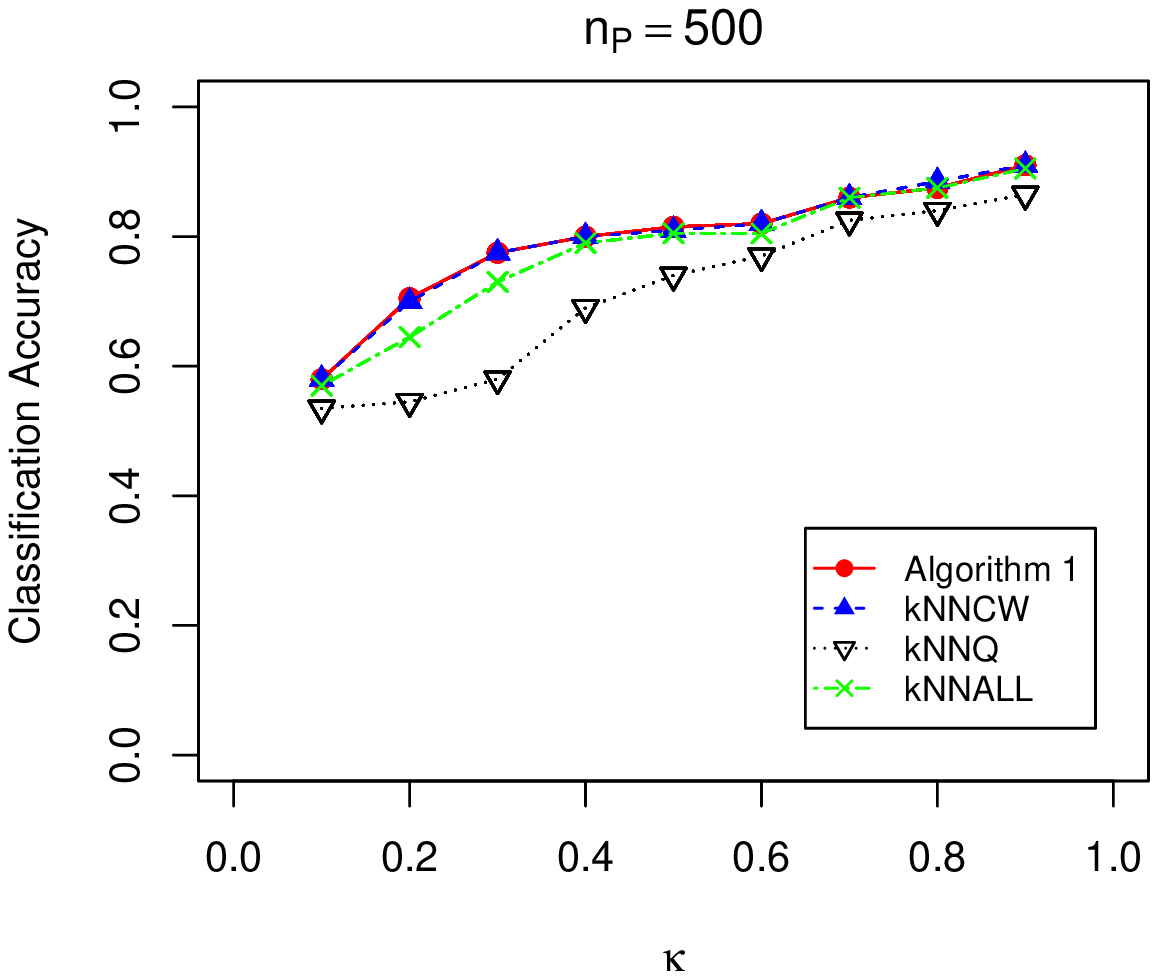}
\includegraphics[width=2 in]{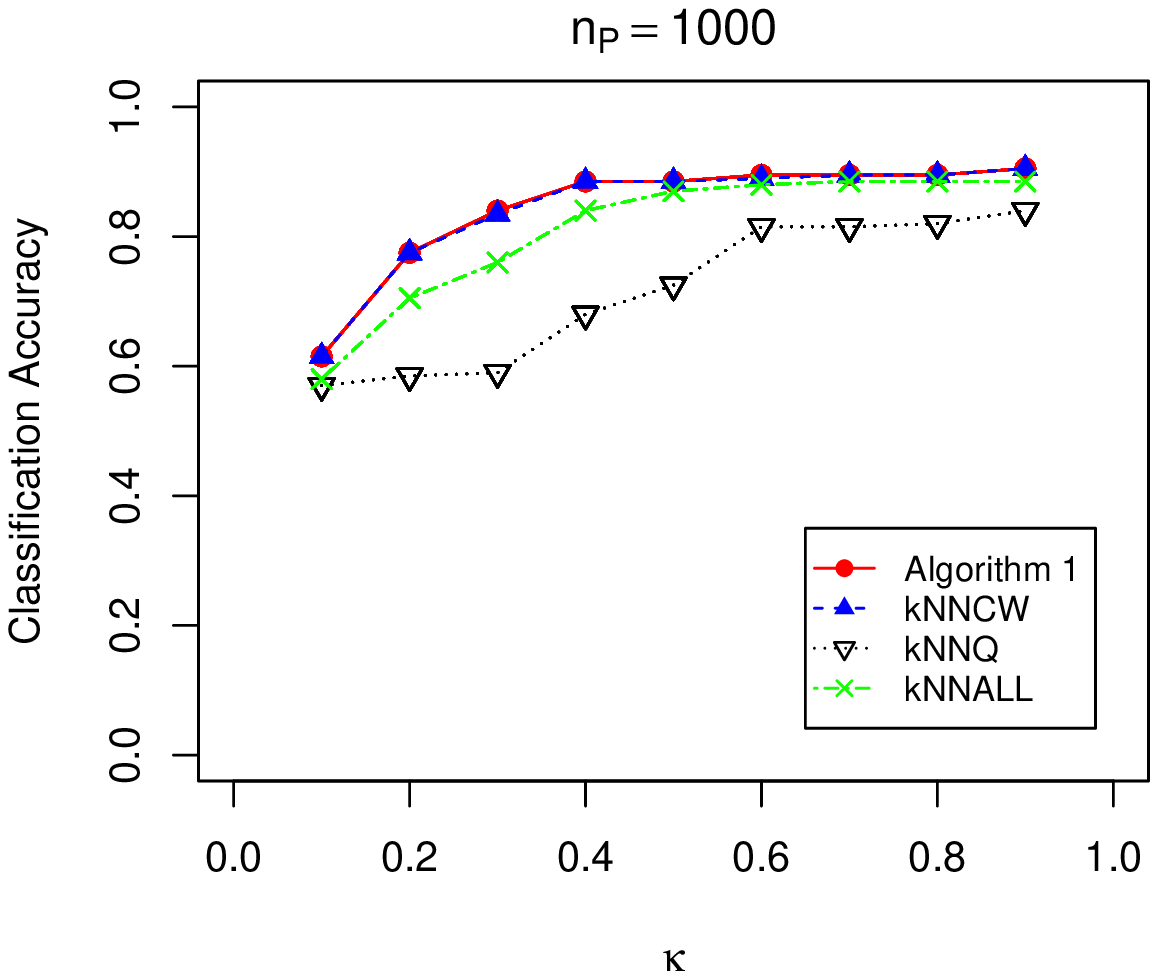}

\includegraphics[width=2 in]{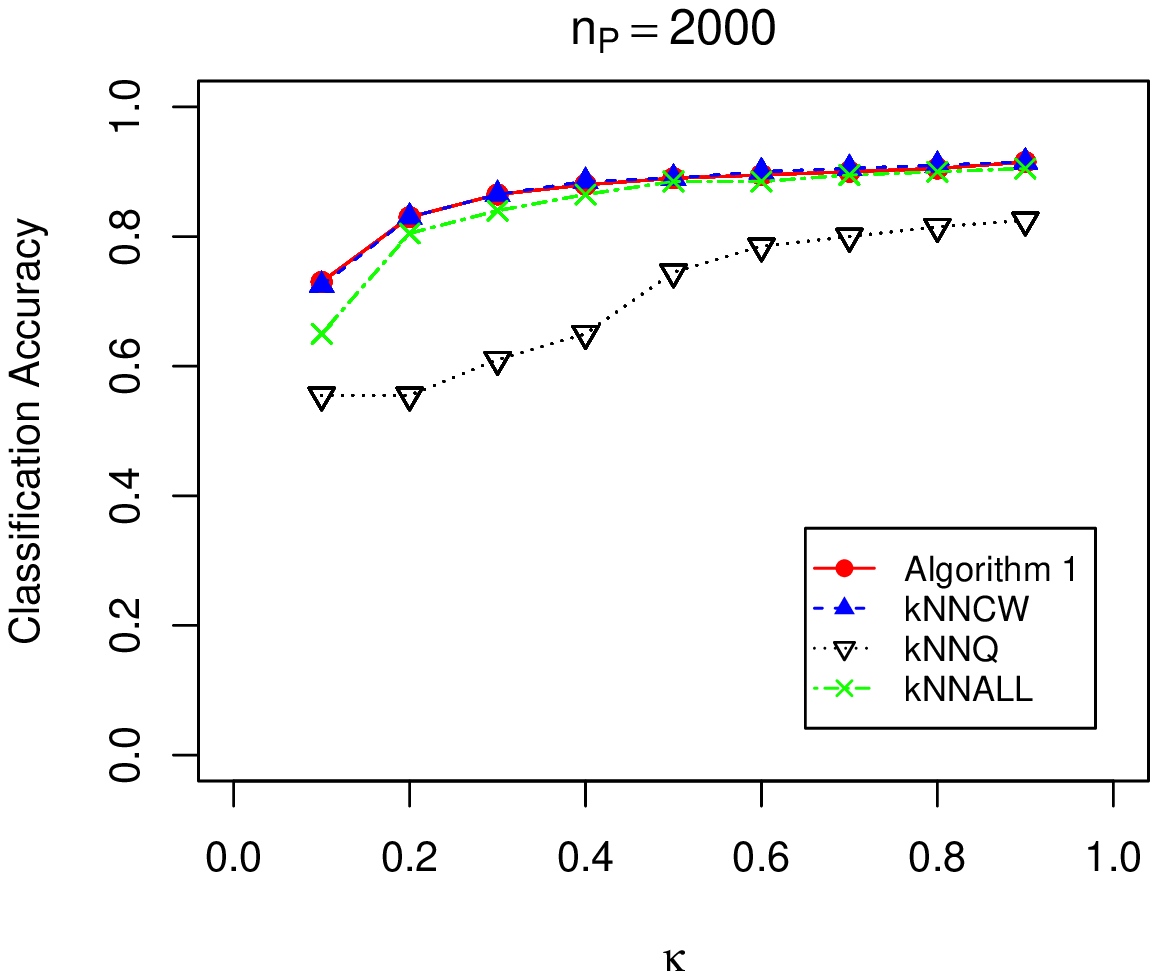}
\includegraphics[width=2 in]{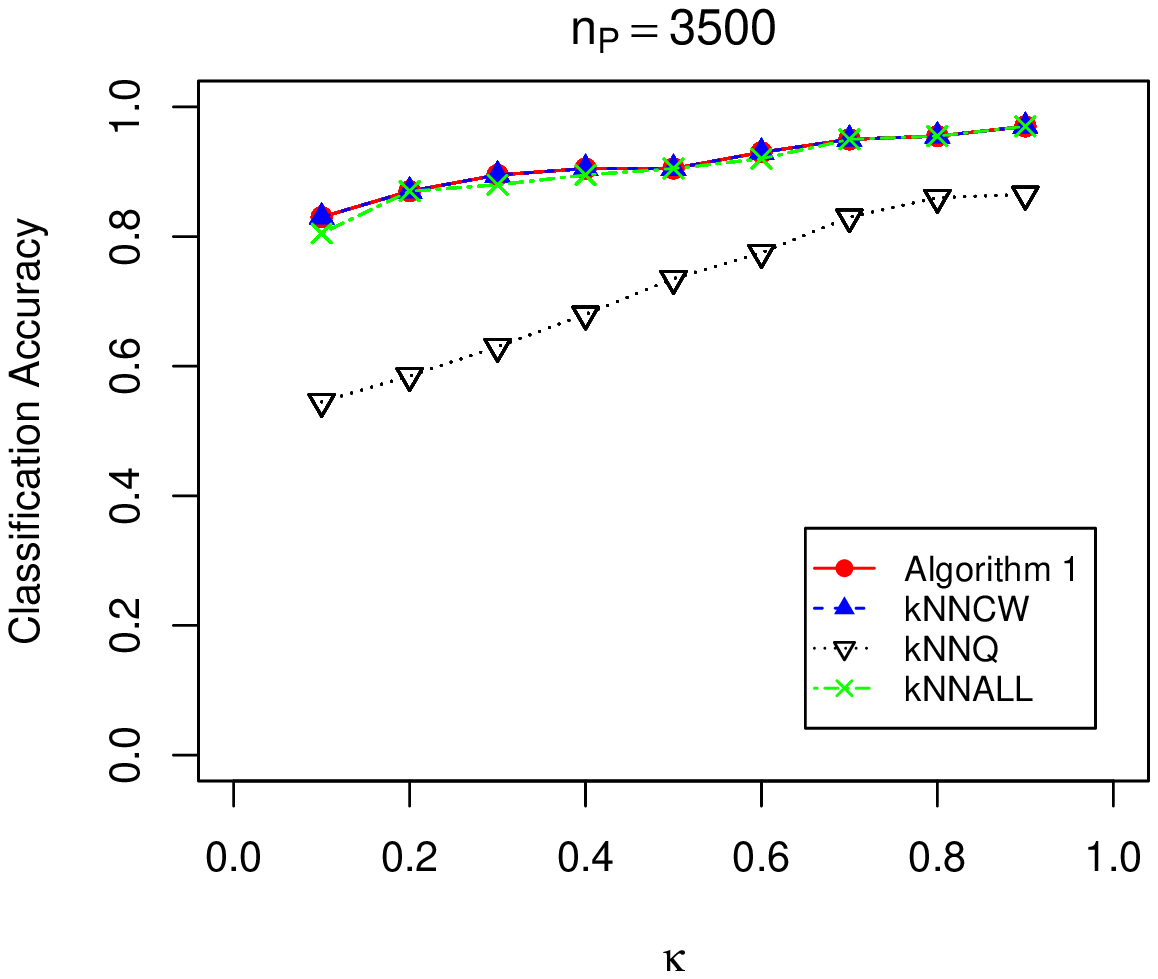}
\includegraphics[width=2 in]{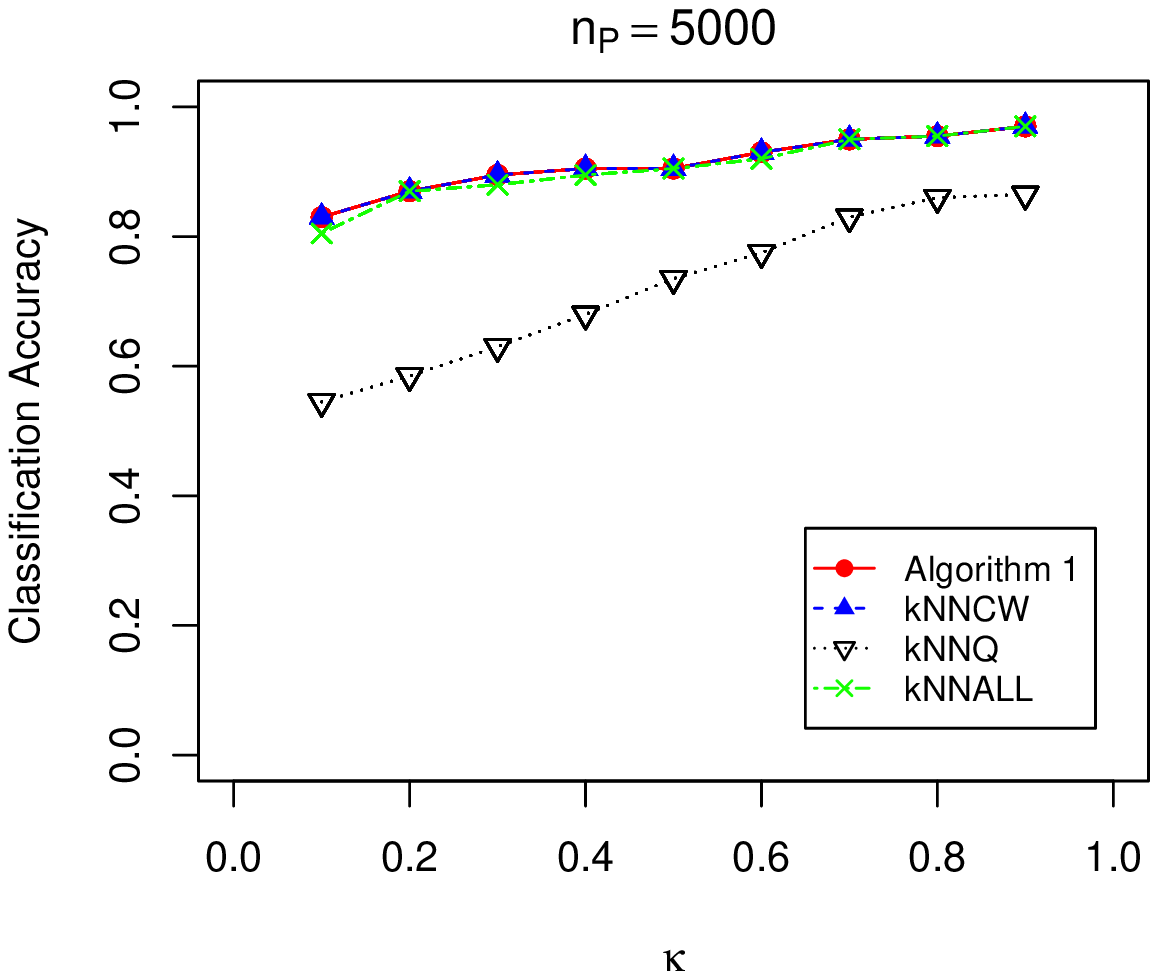}
\caption{\it DGP 2 (Smooth Target): Classification accuracy under different combinations of $(n_P, \kappa)$.}
\label{figure:classification:DGP2}
\end{figure}

\begin{figure}[h]
\centering
\includegraphics[width=2 in]{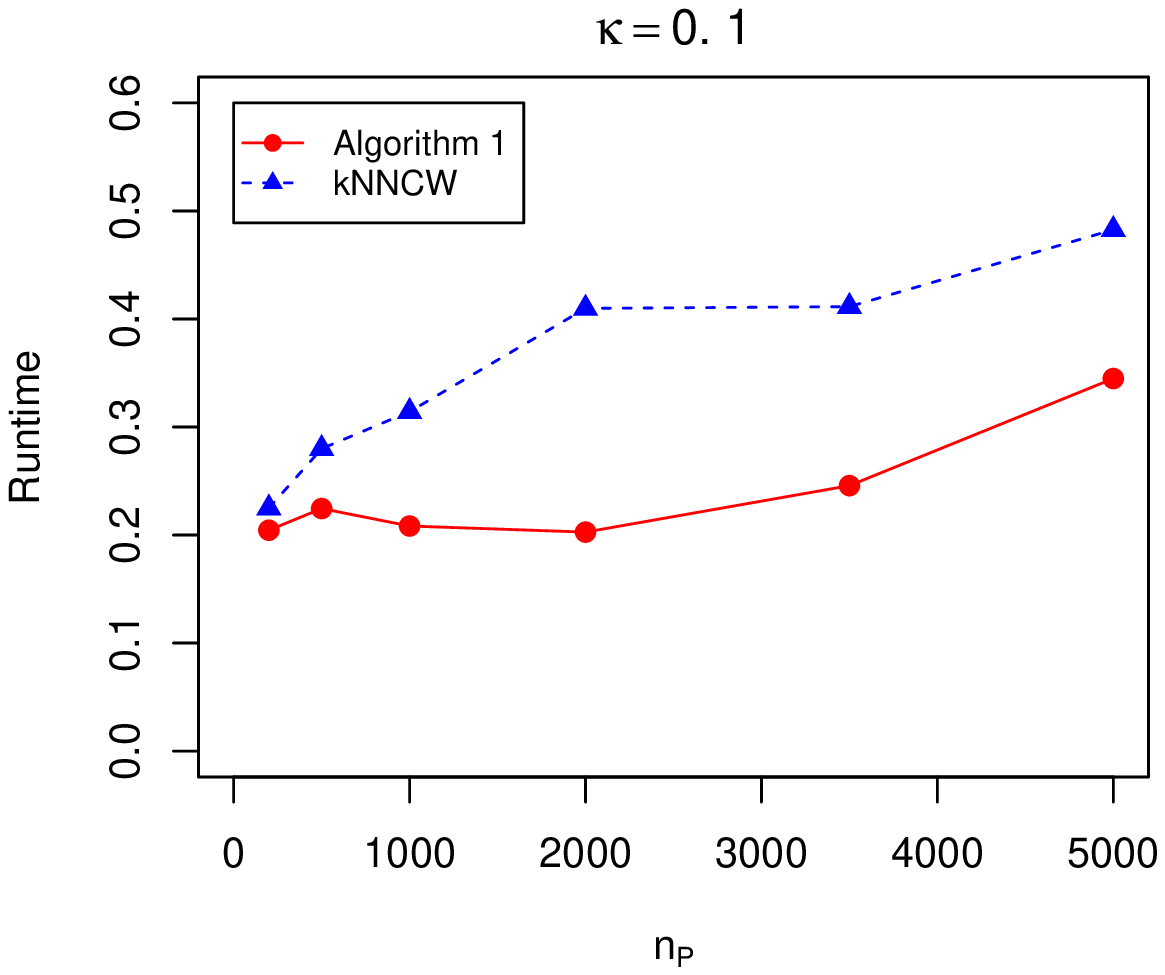}
\includegraphics[width=2 in]{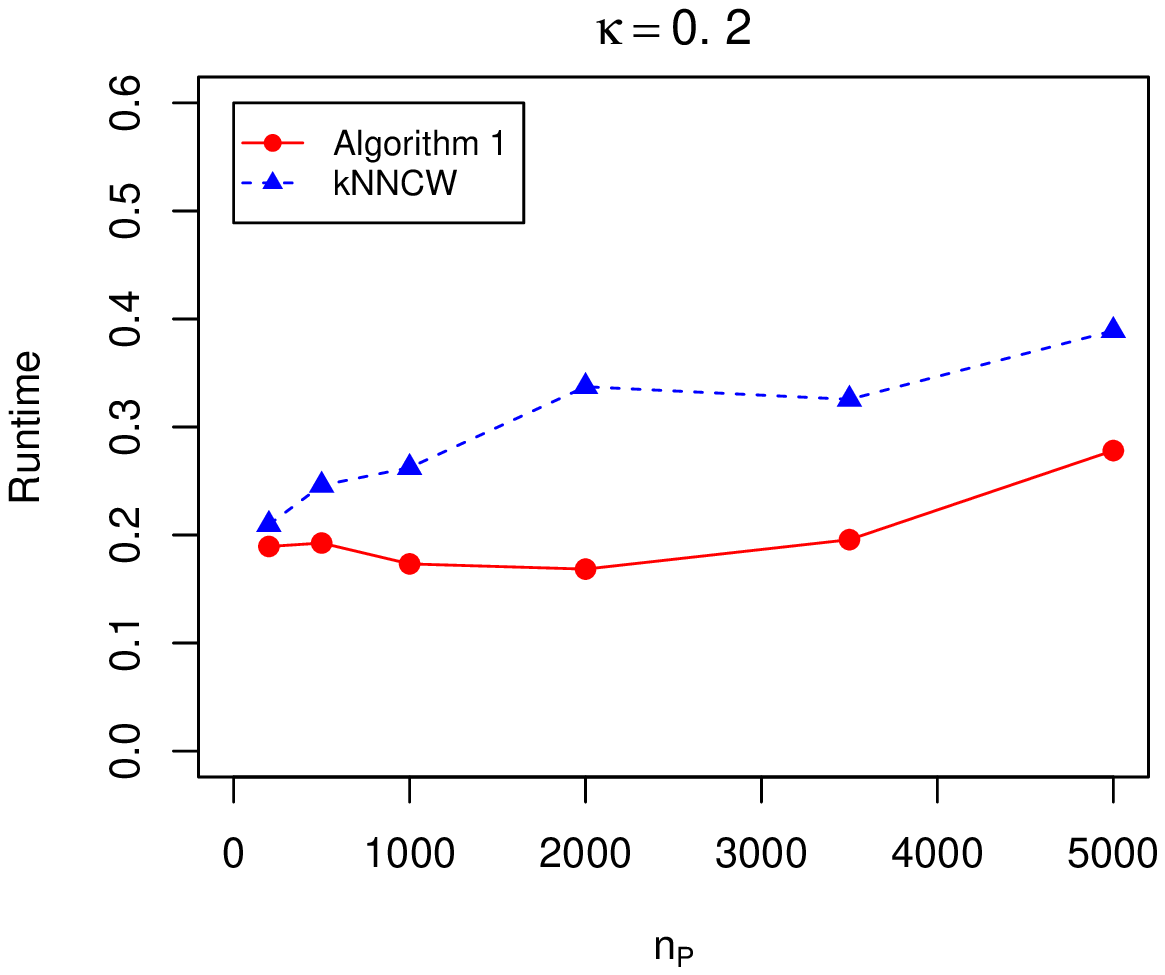}
\includegraphics[width=2 in]{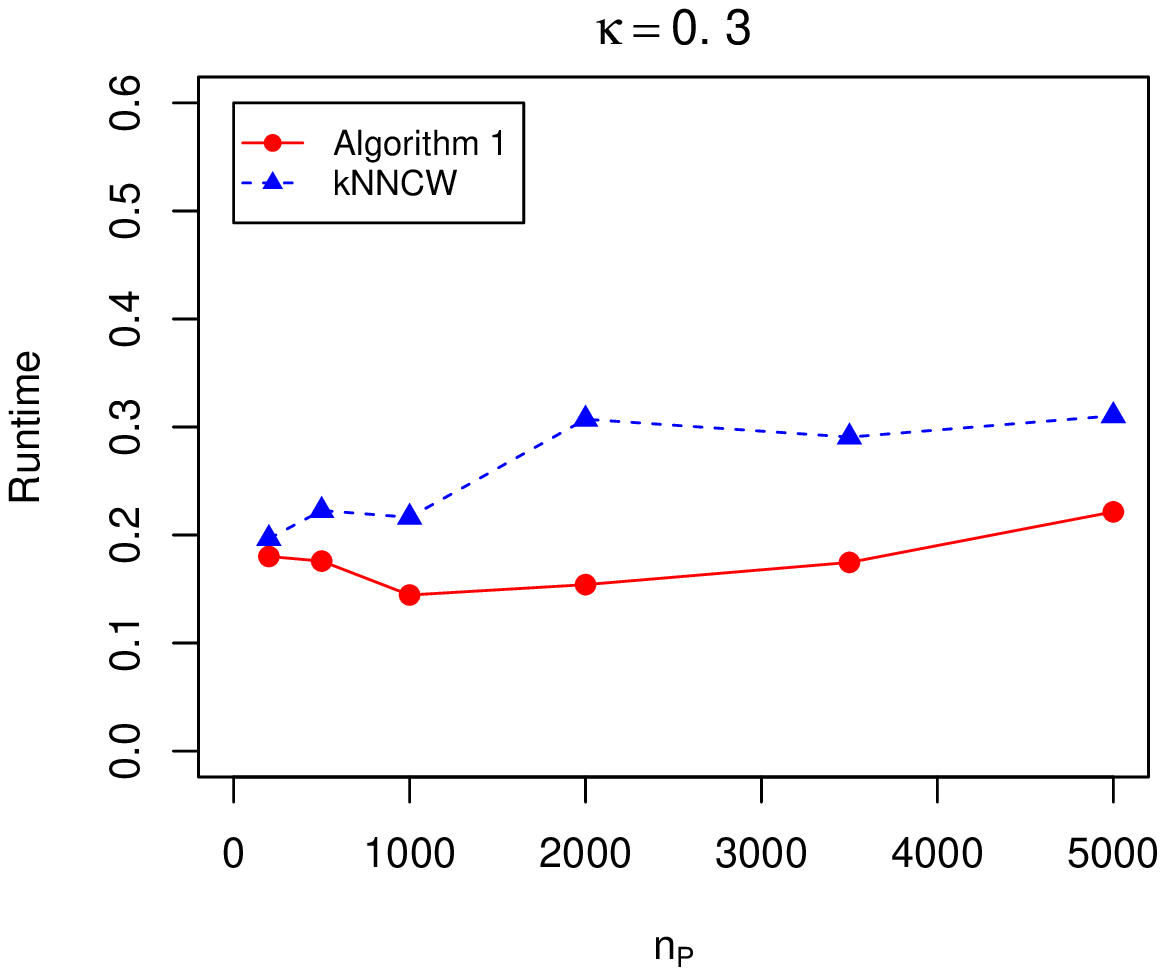}
\includegraphics[width=2 in]{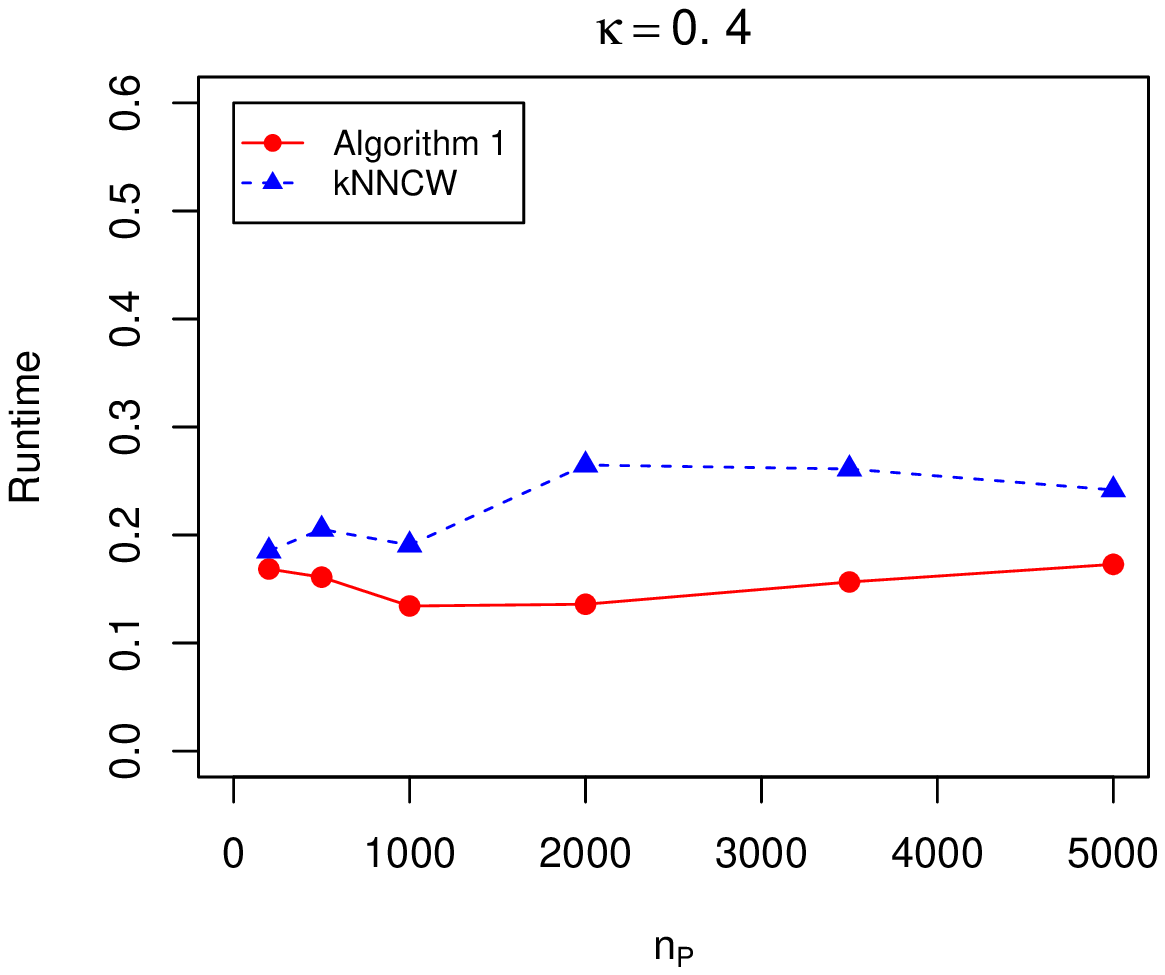}
\includegraphics[width=2 in]{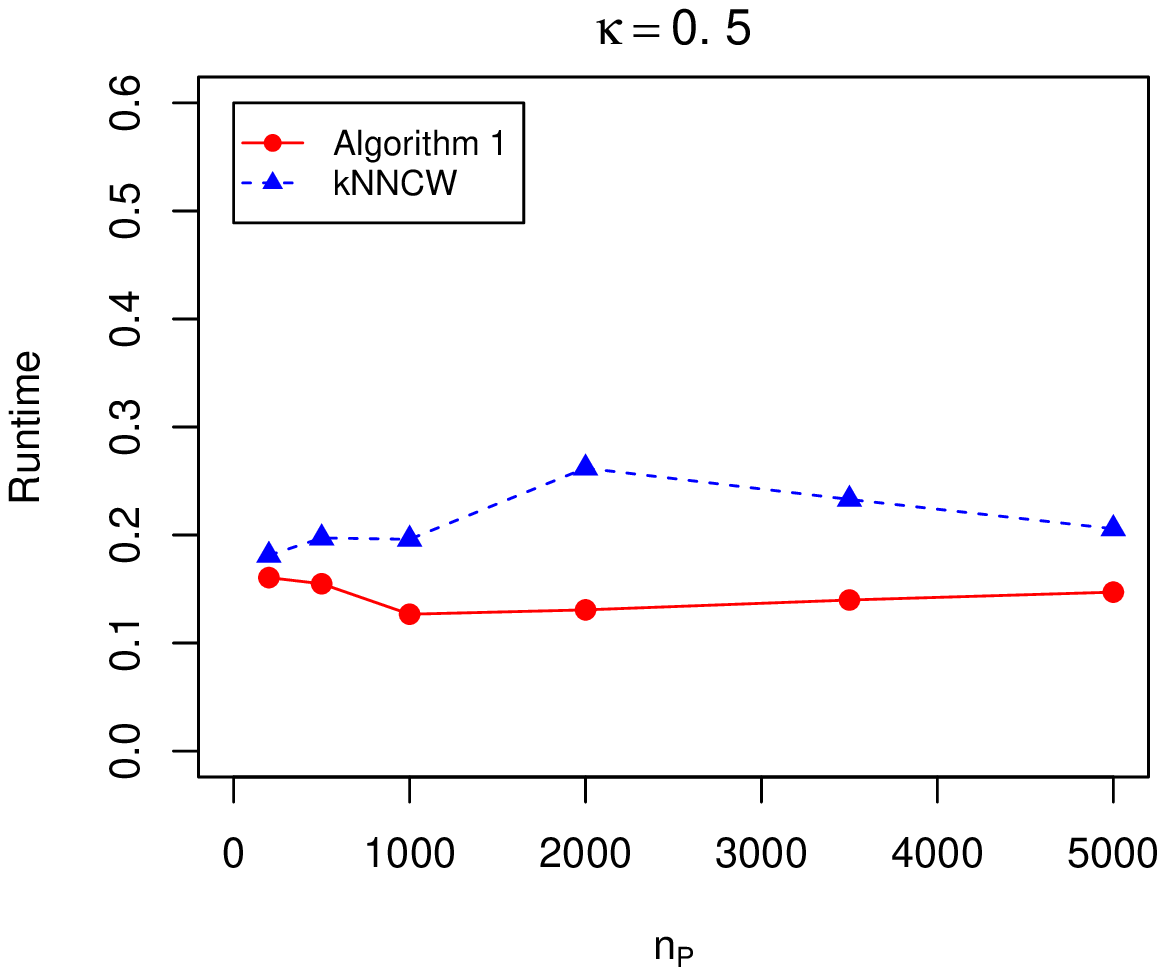}
\includegraphics[width=2 in]{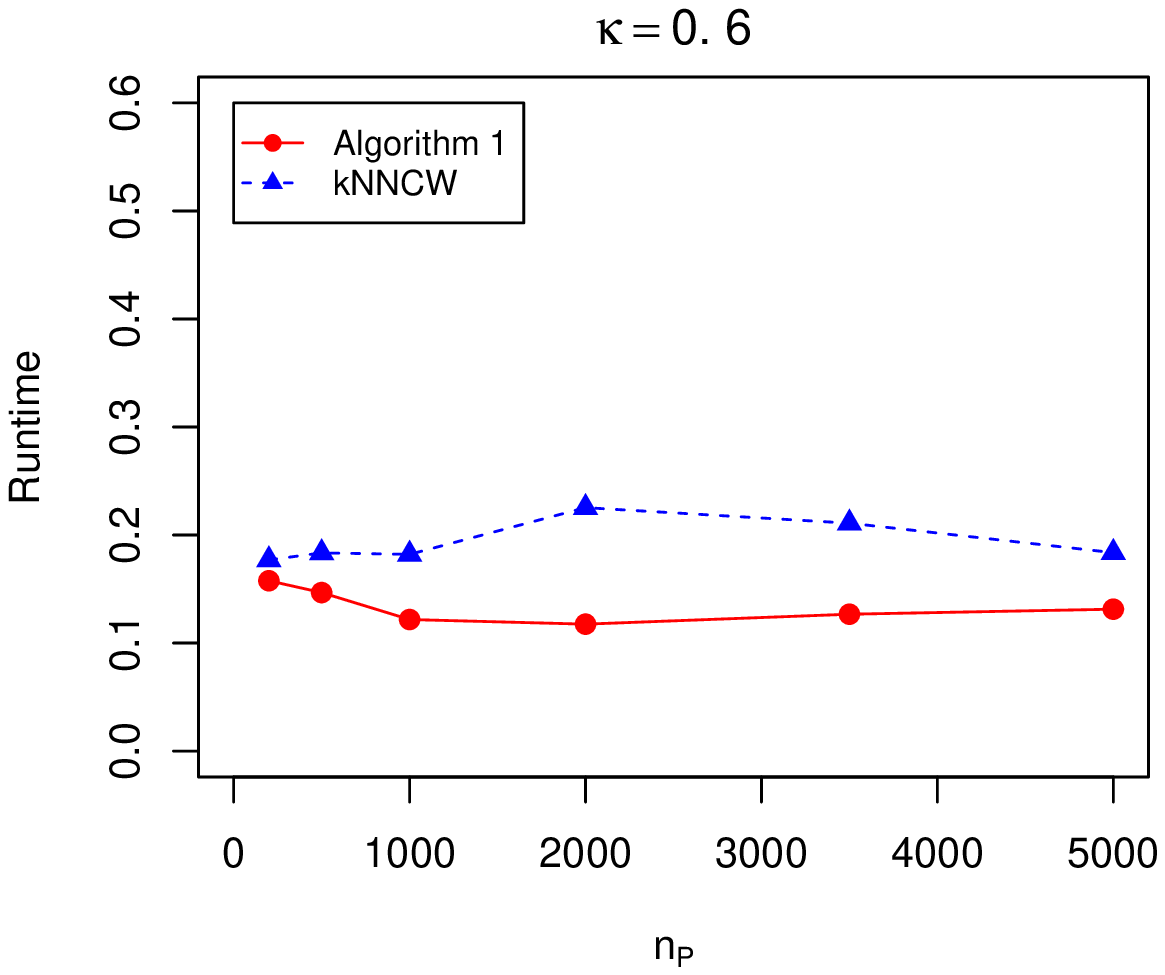}
\includegraphics[width=2 in]{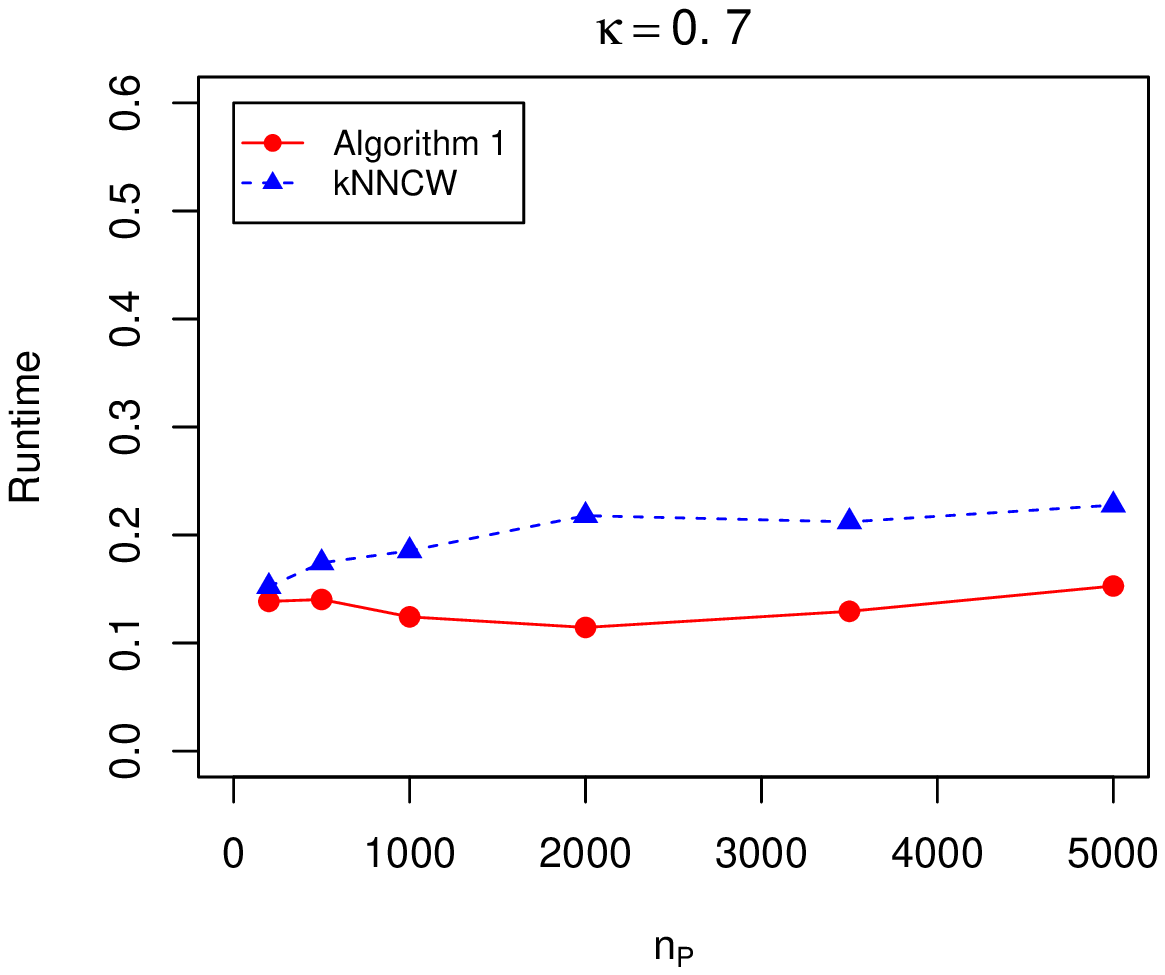}
\includegraphics[width=2 in]{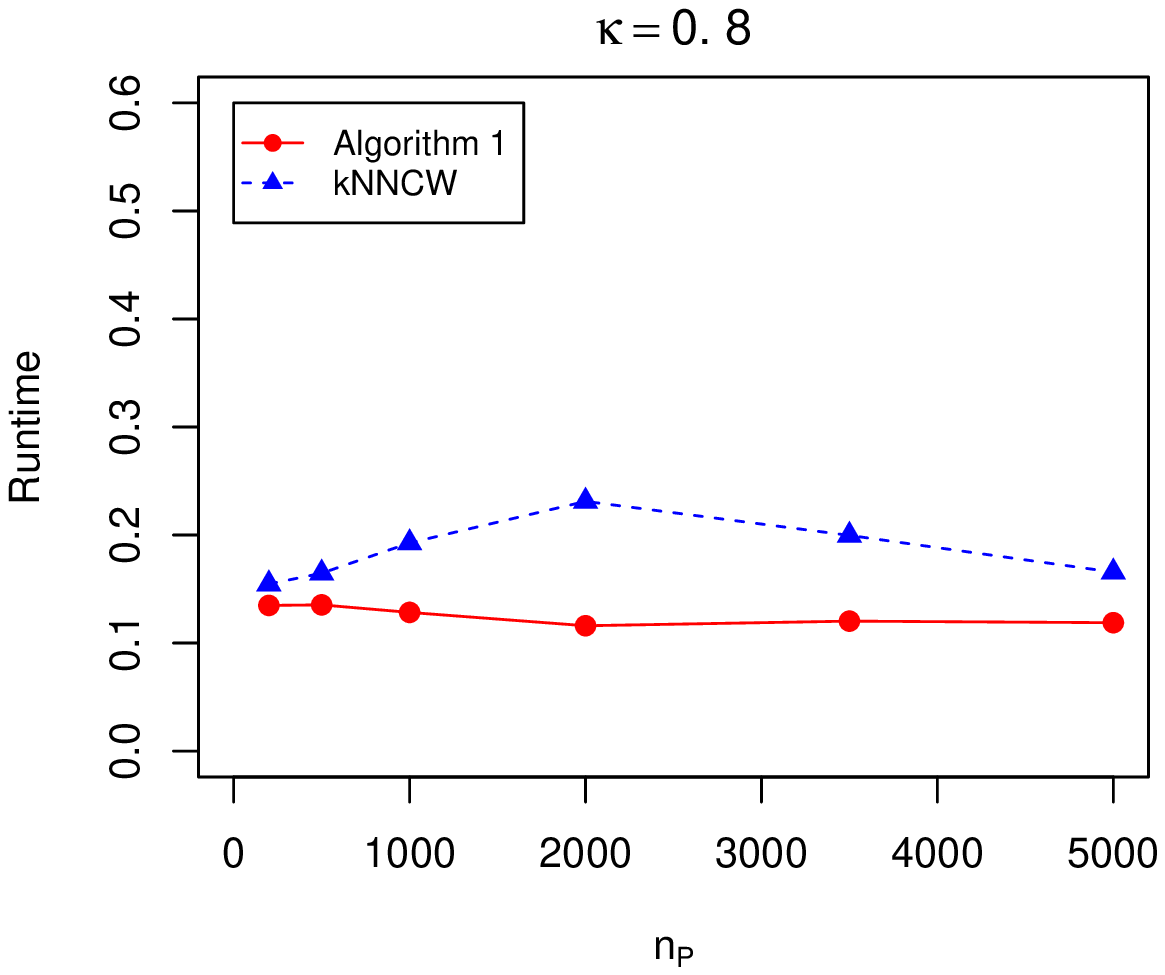}
\includegraphics[width=2 in]{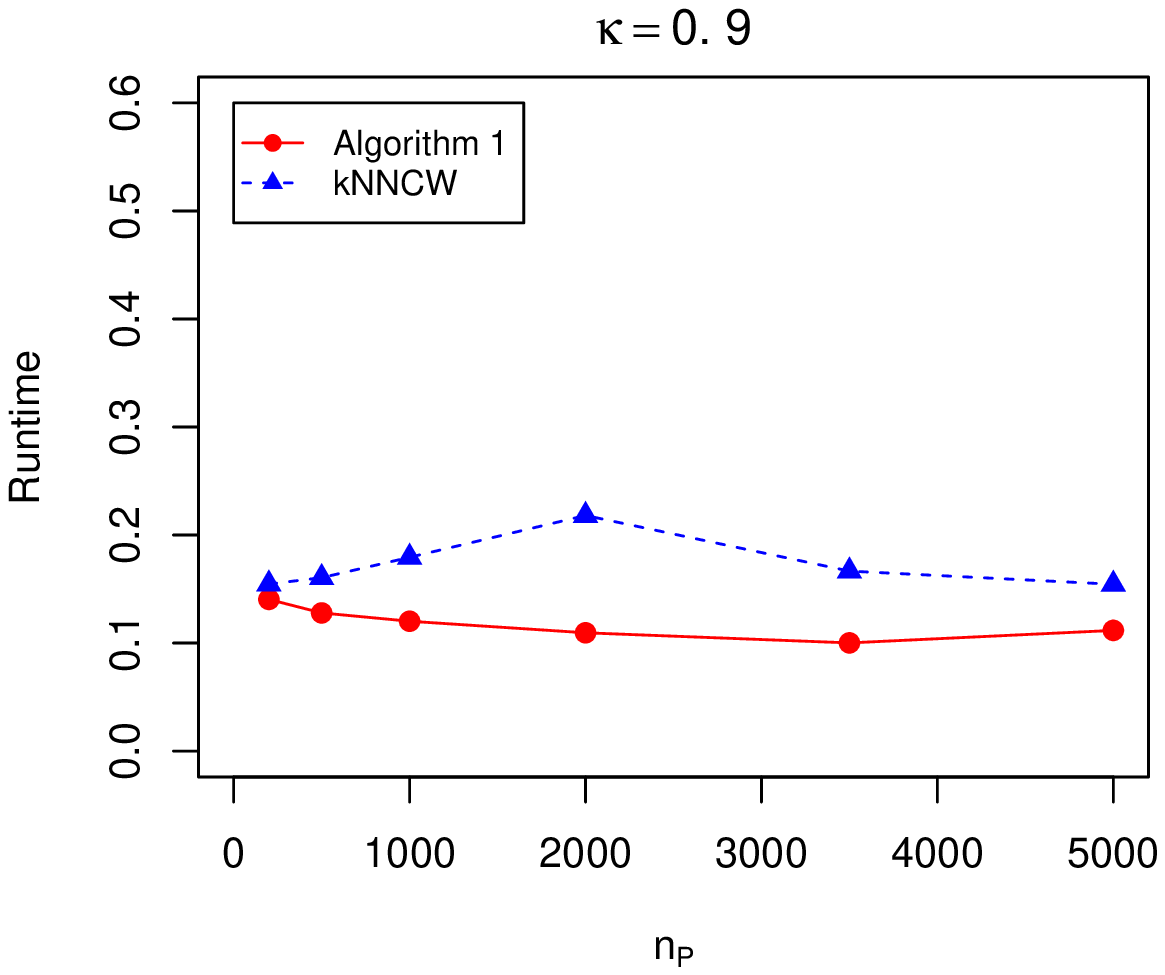}

\caption{\it Runtime (in minutes) of Algorithm \ref{alg:ag2} and kNNCW of DGP 1 under different $(n_P, \kappa)$.}
\label{figure:runtime:DGP1}
\end{figure}

\begin{figure}[h]
\centering
\includegraphics[width=2 in]{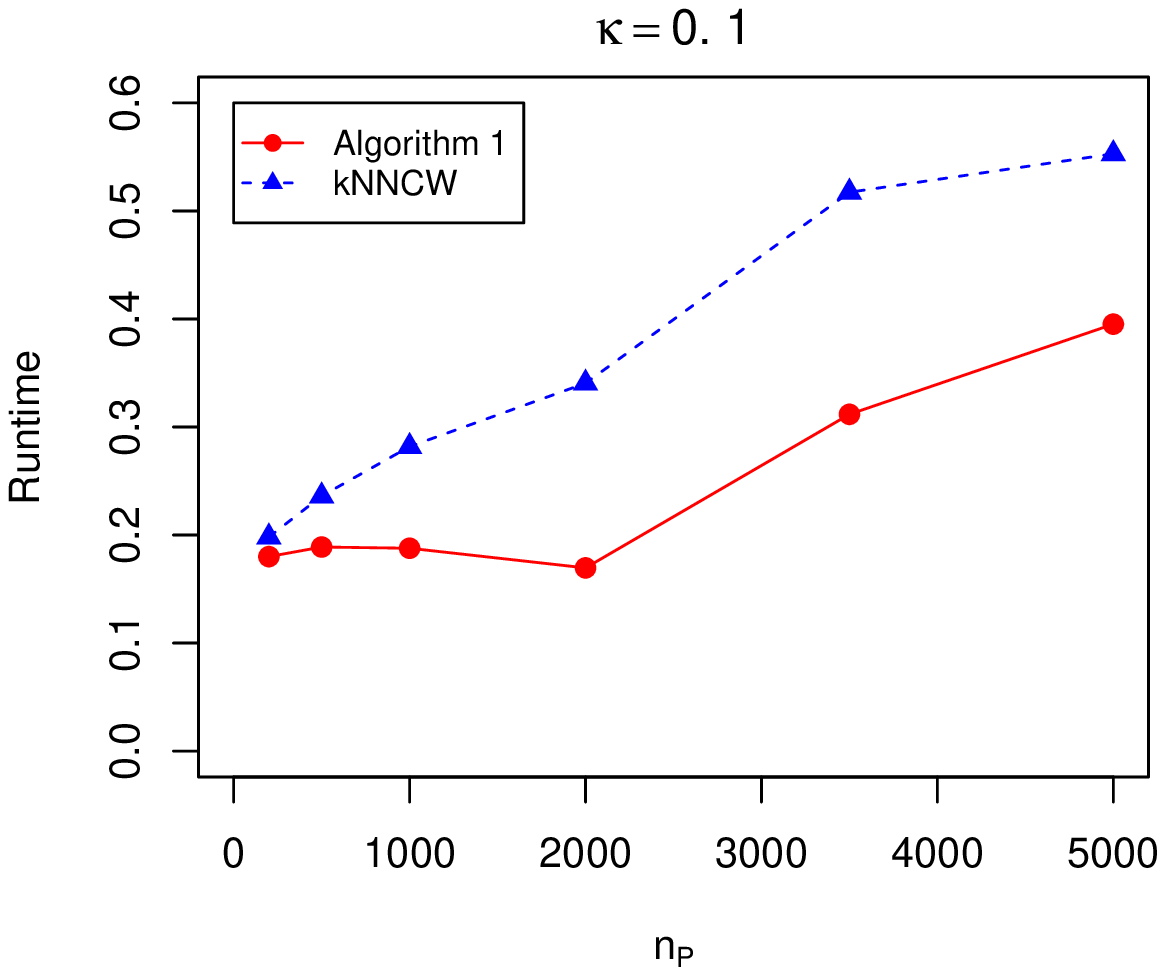}
\includegraphics[width=2 in]{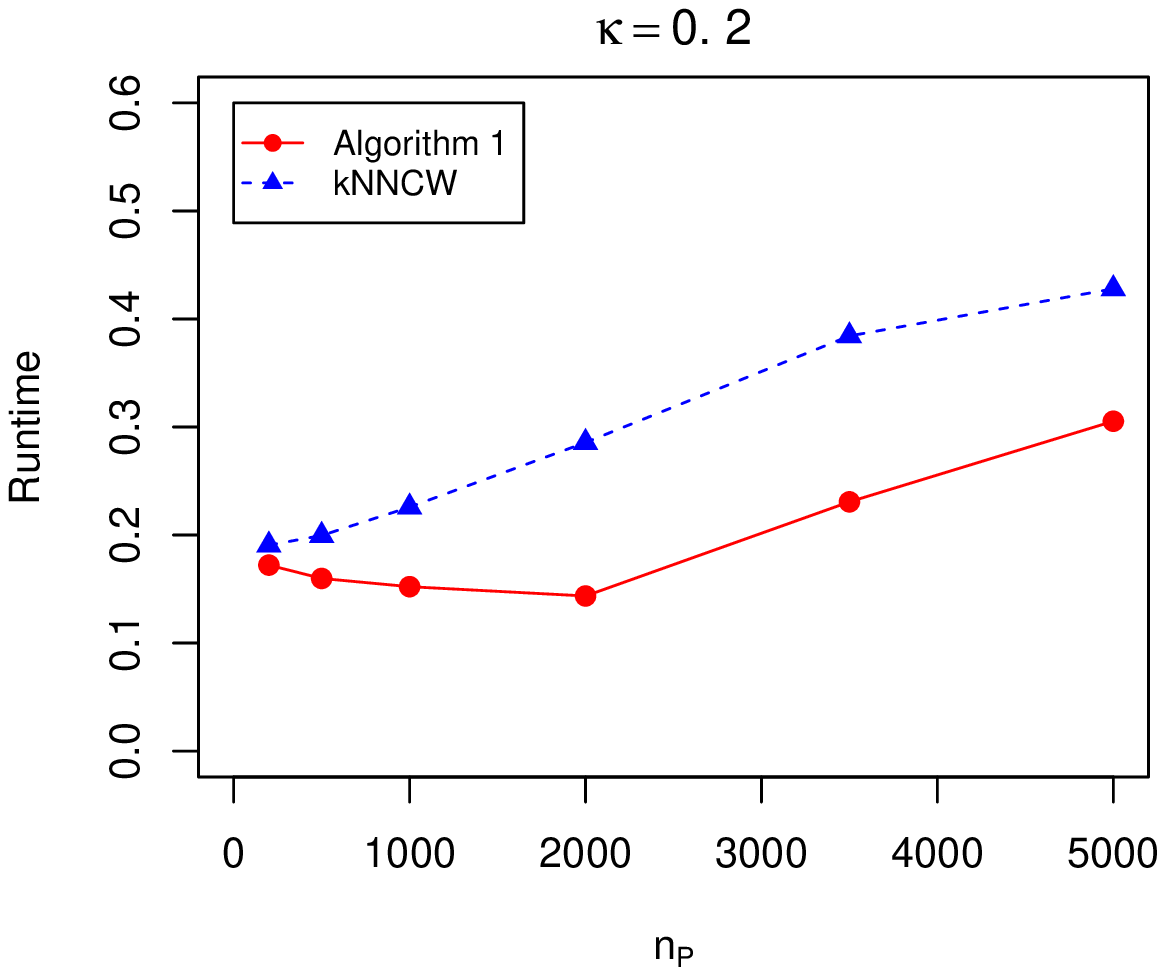}
\includegraphics[width=2 in]{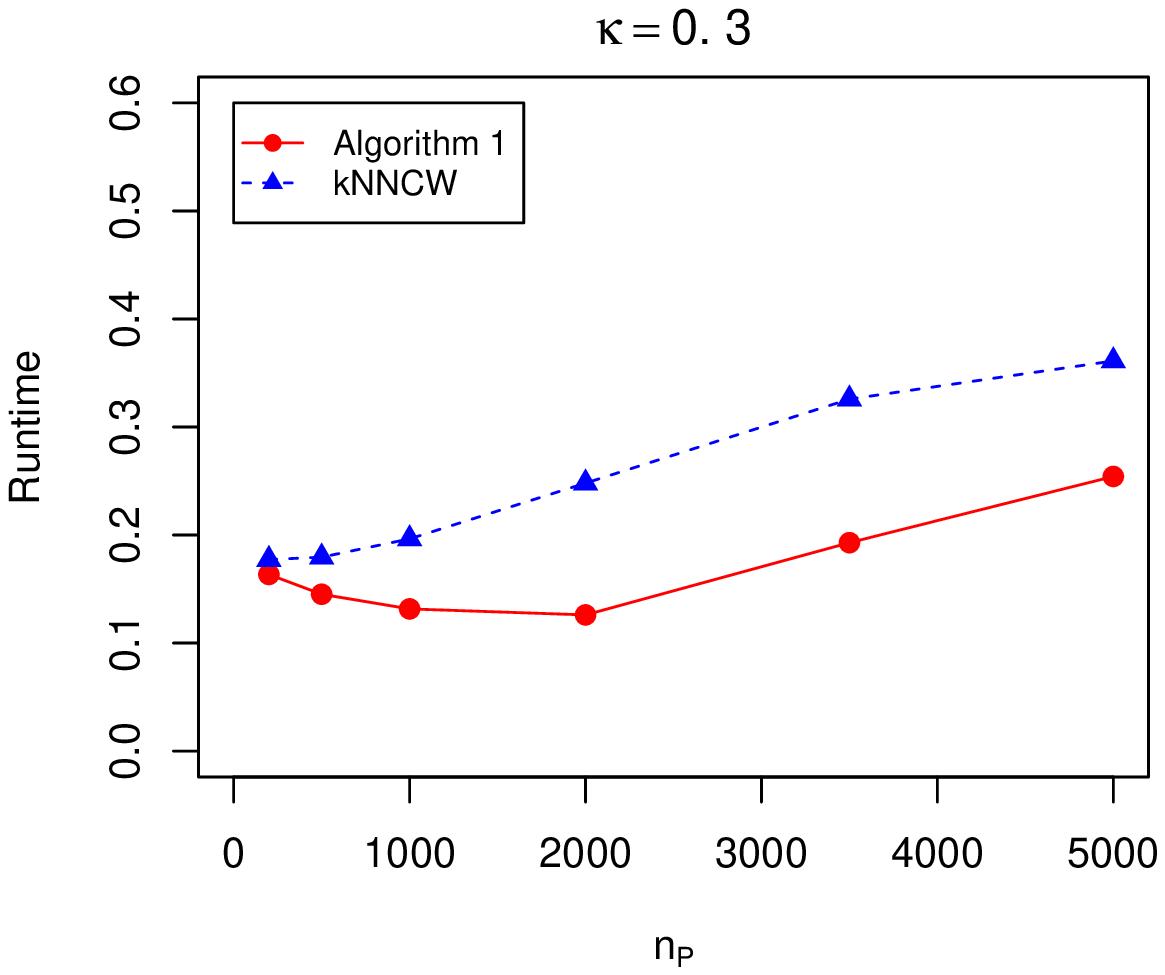}
\includegraphics[width=2 in]{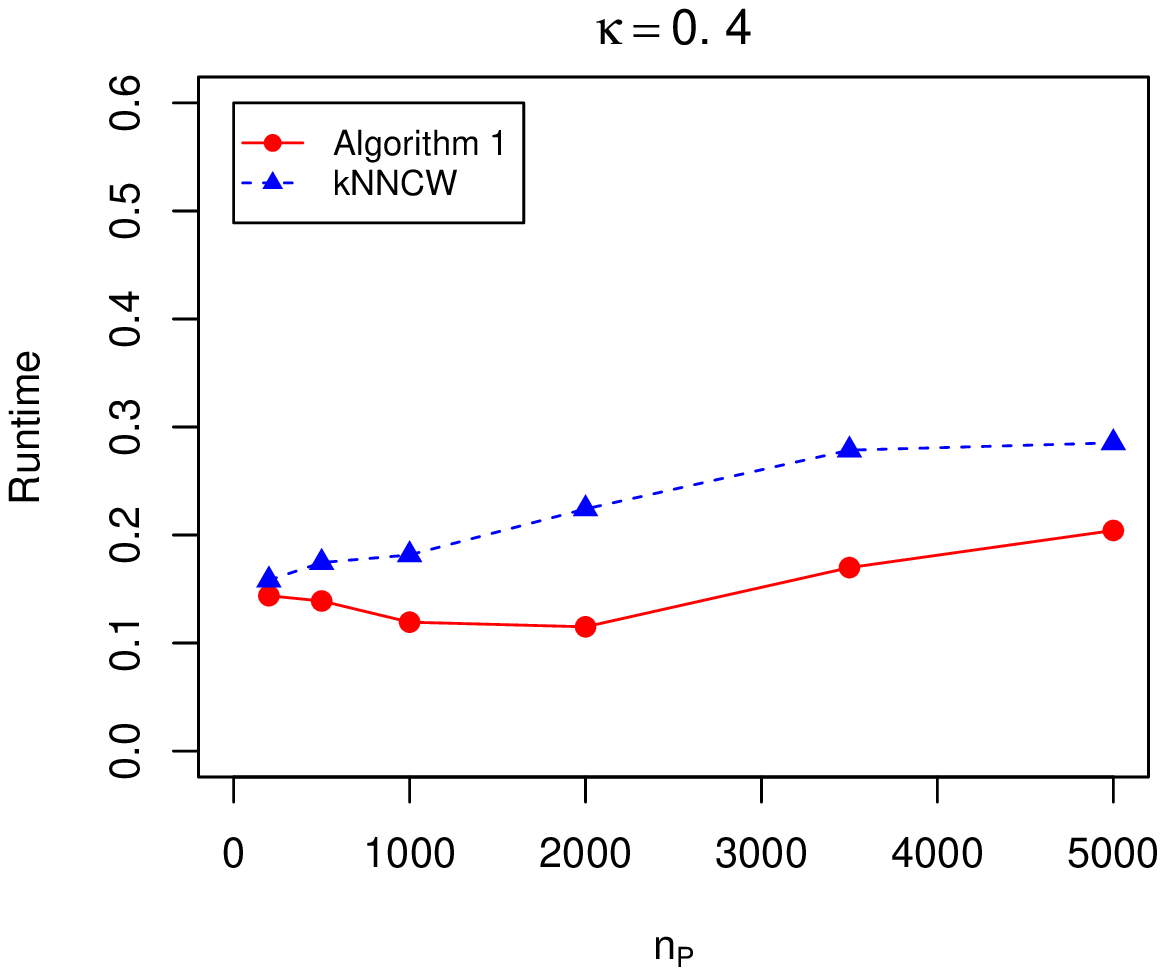}
\includegraphics[width=2 in]{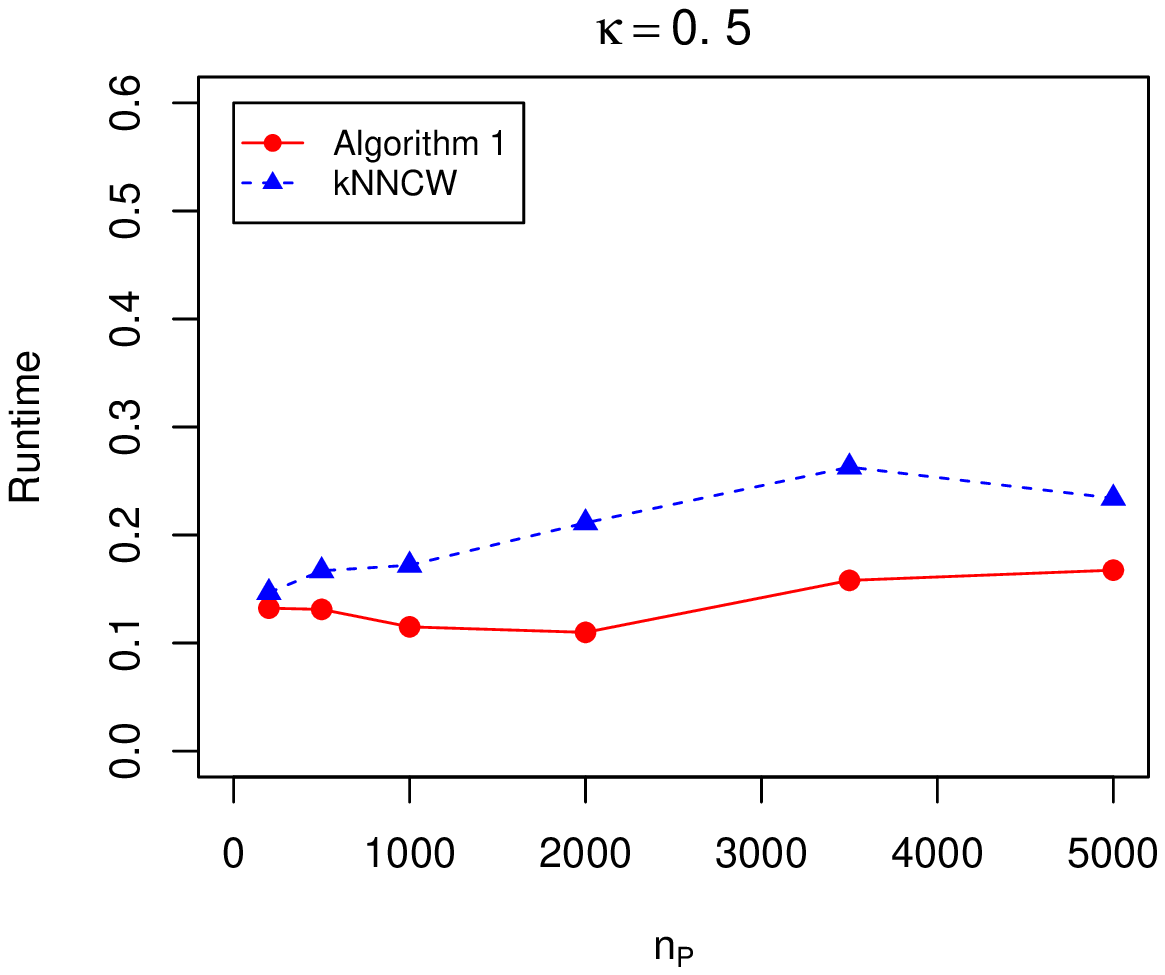}
\includegraphics[width=2 in]{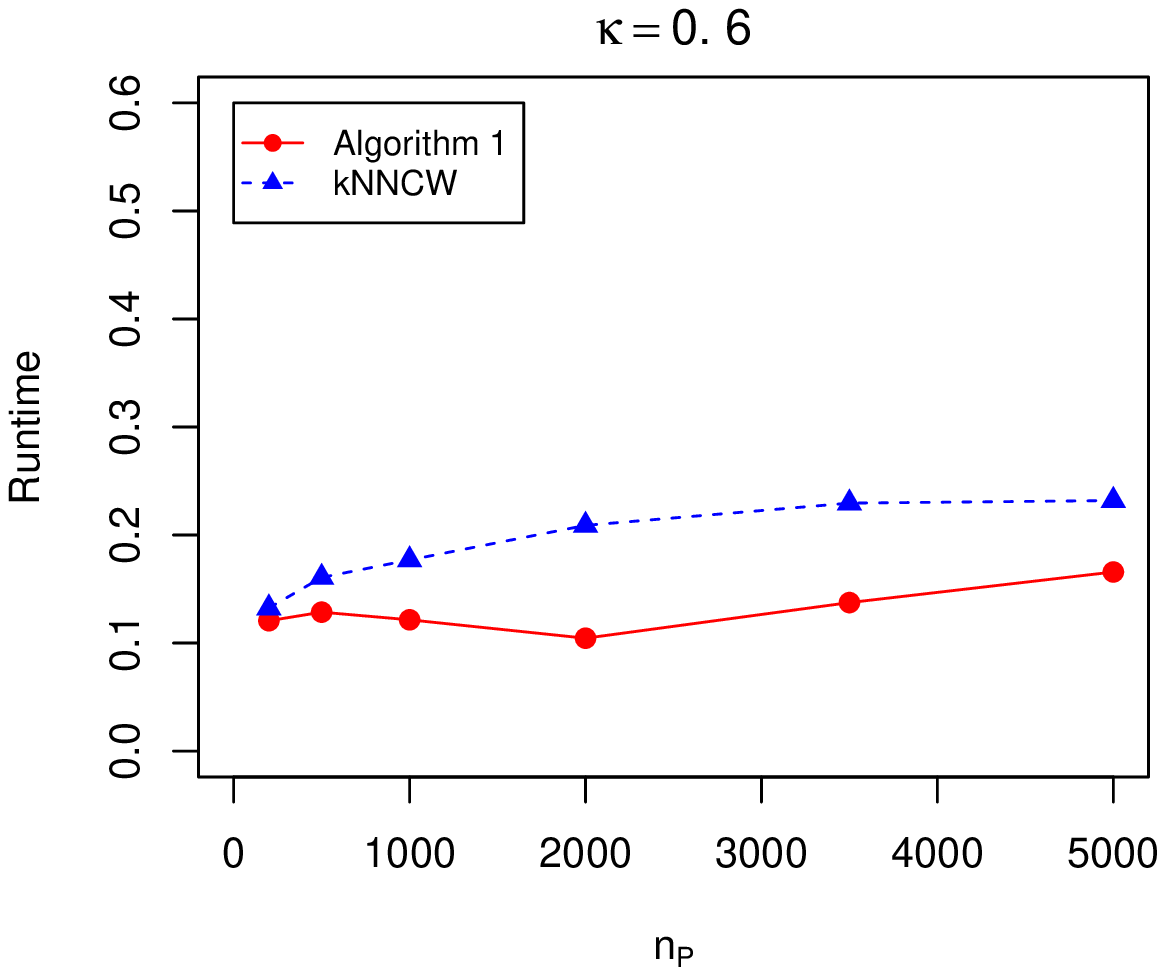}
\includegraphics[width=2 in]{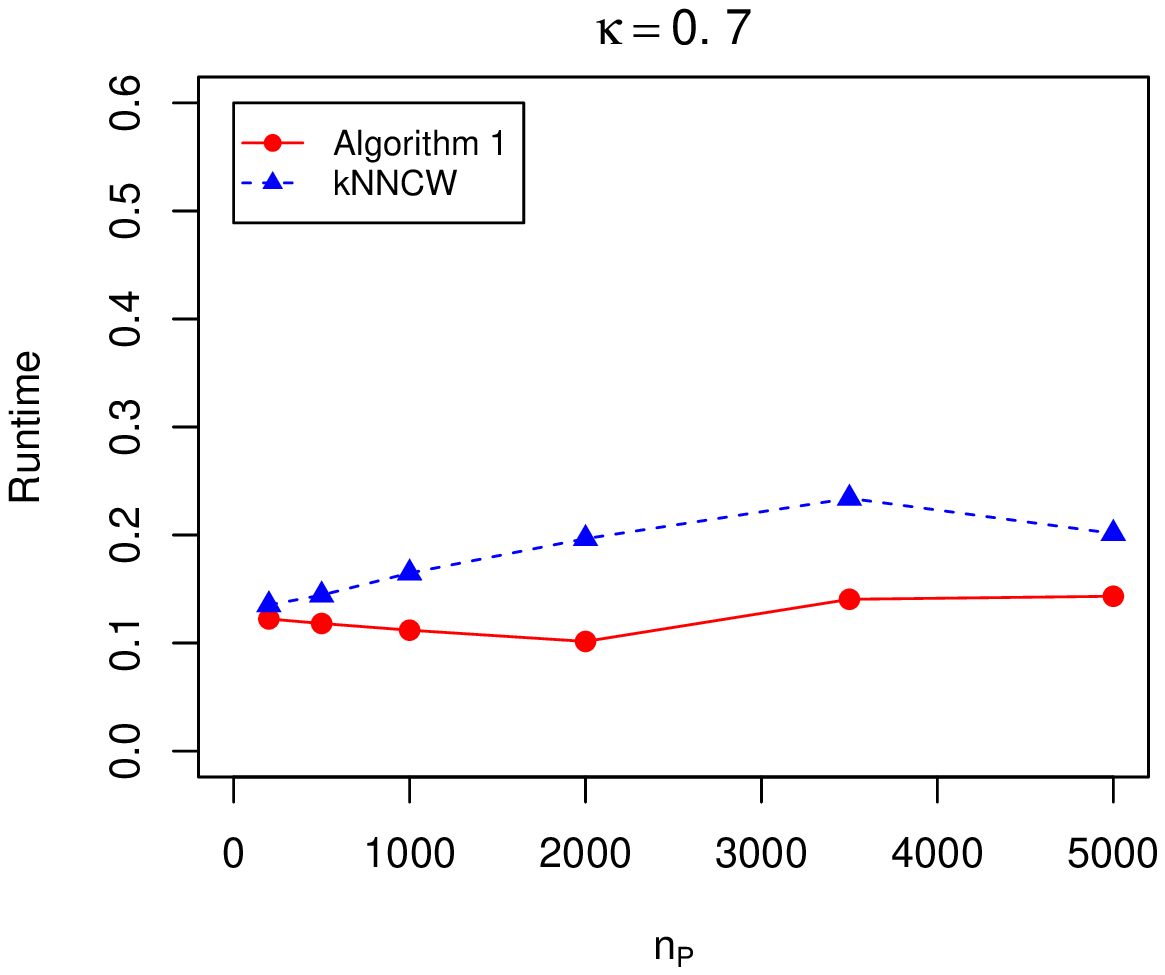}
\includegraphics[width=2 in]{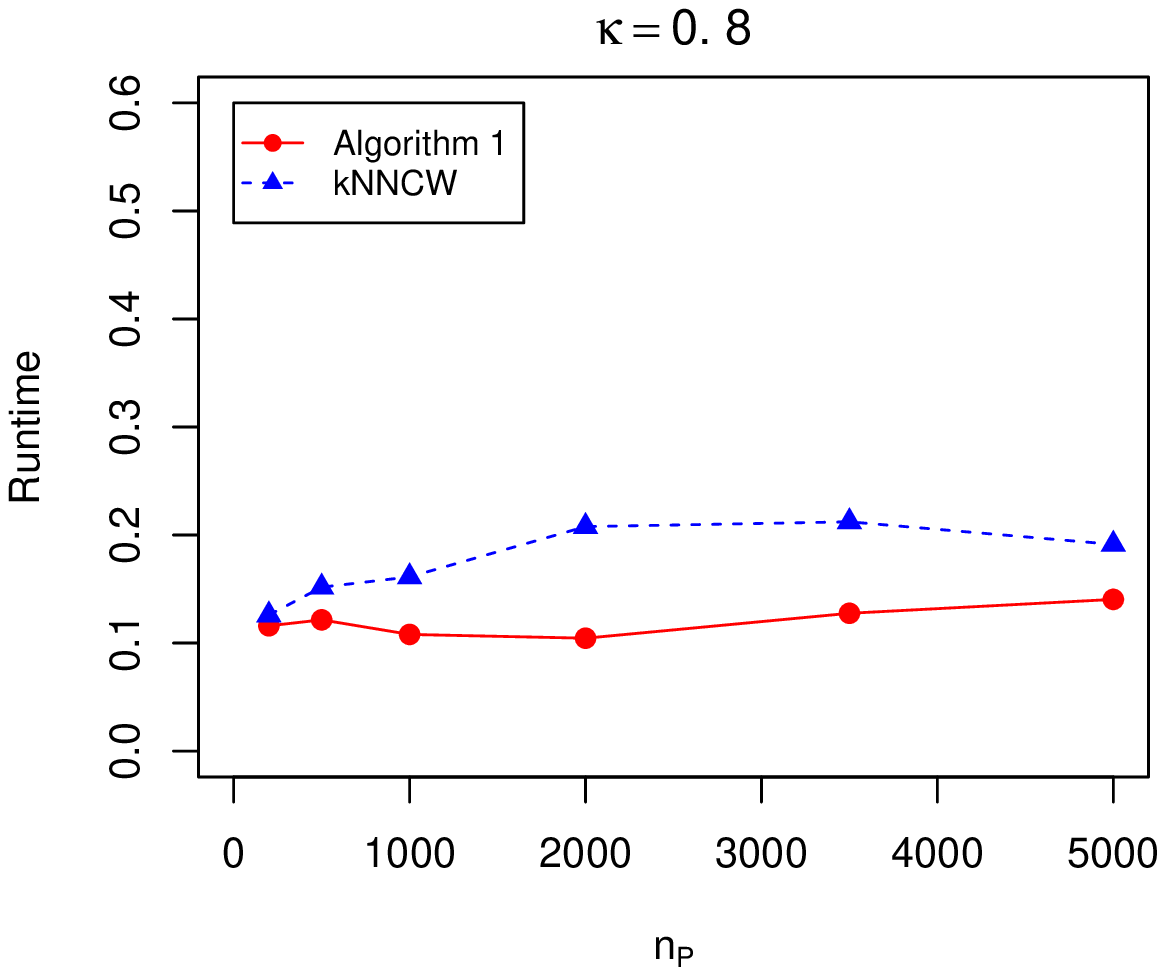}
\includegraphics[width=2 in]{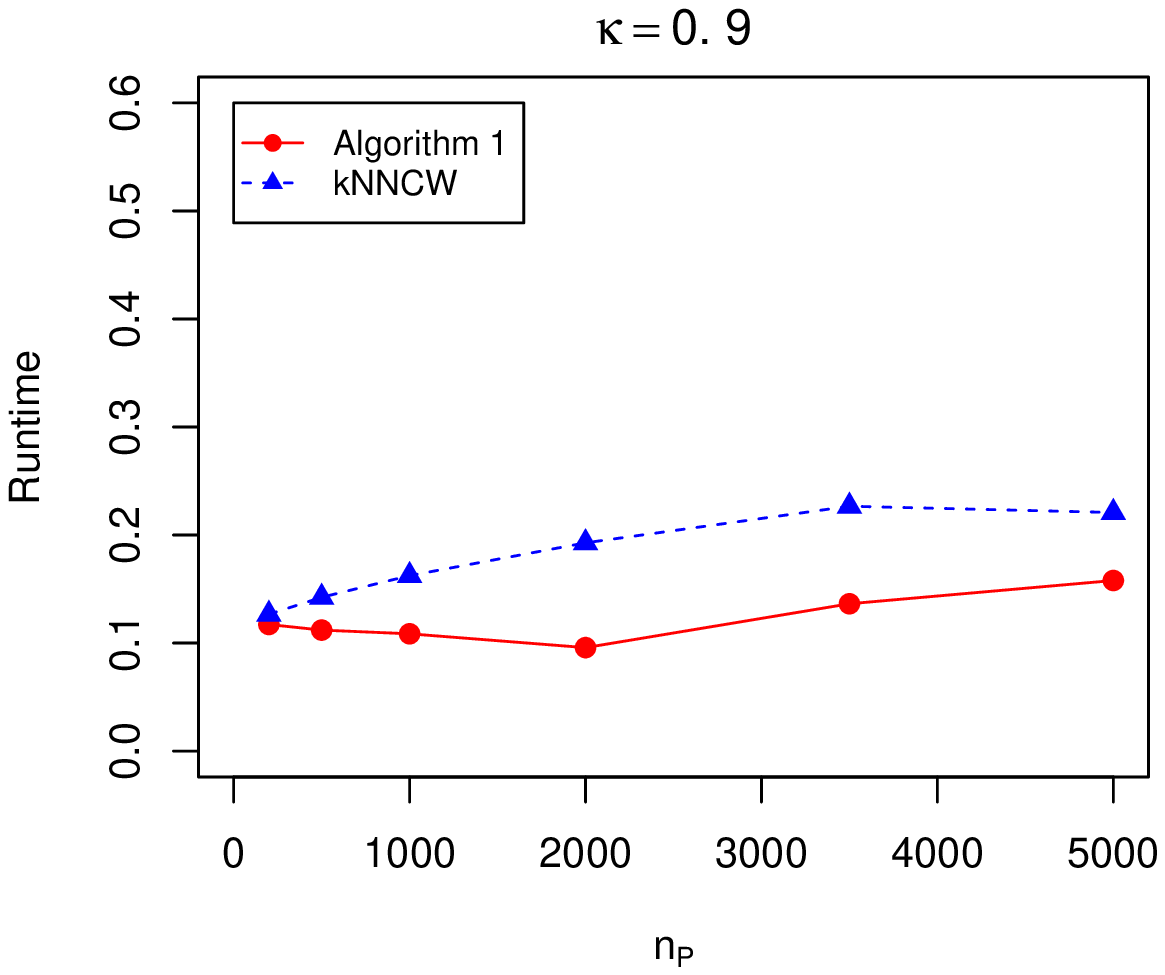}
\caption{\it Runtime (in minutes) of Algorithm \ref{alg:ag2} and kNNCW of DGP 2 under different $(n_P, \kappa)$.}
\label{figure:runtime:DGP2}
\end{figure}

Numerical results are summarized in Figures \ref{figure:classification:DGP1}-\ref{figure:runtime:DGP2}, in which several interesting findings can be observed. First, under different combinations of $n_P$ and $\kappa$, the performance of Algorithm \ref{alg:ag2} is almost identical to $k$NNCW in both DGPs, which meets the theoretical results of Theorem \ref{theorem:two:sample:KNN:regret}. Second, when increasing $n_P$, the performance of $k$NNQ, which is only relying on $Q$-data, becomes much poorer in comparison with the other three classifiers. 

Third, in terms of classification accuracy, $k$NNALL is comparable with Algorithm \ref{alg:ag2} and $k$NNCW  when ether $n_P$ is relatively large or relatively small in comparison with $n_Q$. However, the difference becomes significant when $n_P=1000$ and $n_Q=2000$. Last, according to Figures \ref{figure:runtime:DGP1} and \ref{figure:runtime:DGP2}, our algorithm is faster than $k$NNCW. In particular, the computational advantage of Algorithm \ref{alg:ag2} is notable when $n_P=n_Q=2000$.

\section{Empirical Application}
We apply the proposed adaptive algorithm to the Australian Credit Approval dataset (\citealp{quinlan1987simplifying}) downloaded from UCI machine learning repository (\citealp{Dua:2019}). After removing missing values, we keep four continuous explanatory variables $V_2, V_3, V_7, V_{13}$ and normalize them into $[0, 1]$, whose descriptive statistics are summarized in Table \ref{table:descriptive:statistics}. The response variable $y\in \{0, 1\}$  indicates approval or disapproval status. Based on the binary explanatory variable $V_1\in \{0, 1\}$, we further divide the observations into two datasets: $P$-data consists of $468$ observations and $Q$-data consists of $222$ observations. We randomly selected $n_Q$ observations from $Q$-data with $n_Q=100,120,140$, and combined them with $P$-data to train the four  classifiers: Algorithm \ref{alg:ag2}, $k$NNCW, $k$NNQ, and $k$NNALL. The rest $222-n_Q$ observations are treated as the testing dataset. The classification accuracy is calculated based on $100$ independent replications. Results are summarized in Table \ref{table:credit:classification} which indicate that our proposed algorithm leads to a slightly better classification accuracy.

\begin{table}[H]
\centering
\begin{tabular}{ccccccccc}
\hline\hline
        & $V_2$     & $V_3$     & $V_7$     & $V_{13}$     &  &   & V1  & y   \\ \cline{2-5} \cline{8-9} 
1st Qu. & 0.134 & 0.036 & 0.006 & 0.040 &  & class 0 & 222 & 383 \\
Median  & 0.224 & 0.098 & 0.035 & 0.080 &  & class 1 & 468 & 307 \\
Mean    & 0.268 & 0.170 & 0.078 & 0.092 &  &   &     &     \\
3rd Qu. & 0.360 & 0.257 & 0.092 & 0.136 &  &   &     &    \\\hline
\end{tabular}
\caption{\it Descriptive statistics for the Australian Credit Approval dataset.}
\label{table:descriptive:statistics}
\end{table}

\begin{table}[h]
\centering
\begin{tabular}{ccccc}
\hline\hline
$n_Q$   & Algorithm \ref{alg:ag2}     & KNNCW     & KNNQ     & KNNALL    \\ \hline
%80  & 56.34 & \textbf{56.51} & 51.13 & 56.25 \\
100 & \textbf{57.52} & 56.61 & 52.36 & 56.16 \\
120 & \textbf{57.33} & 56.43 & 52.79 & 56.01 \\
140 & \textbf{56.53} & 56.04 & 53.26 & 55.72\\\hline
\end{tabular}
\caption{\it  Classification performance of the Australian Credit Approval dataset for different classifiers.}
\label{table:credit:classification}
\end{table}

\newpage
\bibliography{ref}{}
\bibliographystyle{apalike}

\end{document}